\documentclass{elsarticle}
\usepackage[
pdfauthor={King},
pdftitle={title},
pdfstartview=XYZ,
bookmarks=true,
colorlinks=true,
linkcolor=blue,
urlcolor=blue,
citecolor=blue,
pdftex,
bookmarks=true,
linktocpage=true,   
hyperindex=true
]{hyperref}

\usepackage{hyperref}
\usepackage{natbib}
\usepackage{amssymb}
\usepackage{amsmath}

\newtheorem{theorem}{Theorem}[section]
\newtheorem{corollary}{Corollary}[section]

\newenvironment{definition}[1][Definition]{\begin{trivlist}
\item[\hskip \labelsep {\bfseries #1}]}{\end{trivlist}}

\def\qed{\hfill$\diamondsuit$}
\newcommand{\ed}{\end{document}}
\newcommand{\bear}{\begin{eqnarray*}}
\newcommand{\eear}{\end{eqnarray*}}
\newcommand{\ben}{\begin{eqnarray}}
\newcommand{\een}{\end{eqnarray}}
\newcommand{\bnn}{\begin{enumerate}}
\newcommand{\enn}{\end{enumerate}}

\newcommand{\beq}{\begin{equation}}
\newcommand{\eeq}{\end{equation}}
\newcommand{\beqs}{\begin{equation*}}
\newcommand{\eeqs}{\end{equation*}}
\newcommand{\bea}{\begin{align}}
\newcommand{\eea}{\end{align}}
\newcommand{\beas}{\begin{align*}}
\newcommand{\eeas}{\end{align*}}
\newcommand{\bcase}{\begin{cases}}
\newcommand{\ecase}{\end{cases}}
\newcommand{\btab}{\begin{tabular}}
\newcommand{\etab}{\end{tabular}}
\newcommand{\bct}{\begin{center}}
\newcommand{\ect}{\end{center}}
\newcommand{\var}{{\rm Var}}
\newcommand{\cov}{{\rm Cov}}

\newcommand{\tr}{{\rm tr}}

\newcommand{\ncip}{\buildrel p \over \nrightarrow}
\newcommand{\cid}{\buildrel d \over \longrightarrow}

\newcommand{\cip}{\buildrel p \over \longrightarrow}
\newcommand{\Ima}{{\rm Im}}

\newcommand{\rank}{{\rm rank}}

\newcommand{\E}{{\rm E}}

\newcommand{\CK}{{\cal K}}
\newcommand{\CA}{{\cal A}}
\newcommand{\CB}{{\cal B}}
\newcommand{\CC}{{\cal C}}
\newcommand{\CD}{{\cal D}}

\newcommand{\CF}{{\cal F}}
\newcommand{\CG}{{\cal G}}

\newcommand{\CM}{{\cal M}}
\newcommand{\CN}{{\cal N}}

\newcommand{\CH}{{\cal H}}

\newcommand{\CO}{{\cal O}}
\newcommand{\CT}{{\cal T}}

\newcommand{\BBC}{{\mathbb C}}
\newcommand{\BBN}{{\mathbb N}}
\newcommand{\BBZ}{{\mathbb Z}}
\newcommand{\BBR}{{\mathbb R}}

\newcommand{\BK}{{\mathbf{K}}}

\newcommand{\BR}{{\mathbf{R}}}

\newcommand{\BX}{{\mathbf{X}}}

\newcommand{\CS}{{\cal S}}

\newcommand{\la}{{\langle}}
\newcommand{\ra}{{\rangle}}

\newcommand{\ba}{\mathbf{a}}
\newcommand{\bb}{\mathbf{b}}

\newcommand{\bbf}{\mathbf{f}}
\newcommand{\bg}{\mathbf{g}}

\newcommand{\boldK}{\mathbf{K}}

\newcommand{\boldX}{\mathbf{X}}

\newcommand{\bolda}{\mathbf{a}}

\newcommand{\boldf}{\mathbf{f}}

\newcommand{\boldg}{\mathbf{g}}

\begin{document}

\begin{frontmatter}
\title{Regularized Functional Canonical Correlation Analysis for Stochastic Processes}

\author{David B. King}
\ead{dbking@indiana.edu}

\address{405 N. Park Ave., Indiana University, Bloomington, IN 47408}

\begin{abstract}
In this paper we derive the asymptotic distributions of two distinct regularized estimators for functional canonical correlation as well as their associated eigenvalues, eigenvectors and projection operators.  The methods we developed utilize regularized estimators which approach the functional operators based in reproducing kernel Hilbert spaces (RKHS) as the regularization parameter approaches zero. In addition to providing some justification for the RKHS methods, we explore the asymptotics of regularized operators associated with both Tikhinov and truncated singular value decomposition (TSVD) type regularization. Together, these regularization methods represent two of the most commonly utilized forms of regularization.
\end{abstract}

\begin{keyword}
Canonical Correlation; Asymptotic Distributions; Stochastic Processes; Reproducing Kernel Hilbert Spaces; Regularization; Inverse Problems  \\

{\it AMS 2000 Subject Classification:} Primary 62H20 \sep 60E05 \sep 62H25 \sep 62M99 \sep 45B05 \sep 45Q05
\end{keyword}

\end{frontmatter}

\section{Introduction}

The goal of multivariate canonical correlation analysis (MCCA; Hotelling \citep{Hotelling36}) \nocite{Wichern02} \nocite{Anderson:03} is to identify and quantify the associations between two random vectors $\boldX_1 \in \BR^{N_1}$ and $\boldX_2 \in \BR^{N_2}$. Recently interest has been focused on the extension of this notion to the collection and analysis of ``functional data'' where the term refers to observations that are curves or sample paths of continuous time stochastic processes.  Although development of statistical methodology for the analysis of functional data has been an active research area for well over twenty years, the current popularity of functional data analysis (FDA) is due, in large part, to monographs by Ramsay and Silverman \citep{Ramsay02} \citep{Ramsay05}. What separates functional data from ordinary multivariate data is that the observed data are sample paths from stochastic processes $X_1(\cdot)$ and $X_2(\cdot)$, which are assumed to be elements of some infinite dimensional and separable Hilbert space consisting of functions defined on an index set $E$, such as $[0,1]$ or $\BBZ$.  In this setting, the covariance matrices which are central to the development of the theory of MCCA are replaced by covariance operators of integral type.  In this infinite dimensional case, some difficulties regarding the definition of the sample canonical correlation have already been observed in Leurgans et al. \citep{lms93}. These authors argue that some kind of smoothing or regularization is indispensable when dealing with estimating the sample canonical correlation. The source of the difficulty in the functional data case is that the sample estimators for covariance operators have finite rank while, in principal, they operate in an infinite dimensional Hilbert space. Leurgans et al. \citep{lms93} points out that as a consequence, the sample principal canonical correlation will always be 1 if no regularization or smoothing is done. This problem originates from the fact that when the number of time points at which the processes are measured becomes larger than the sample size, it will always be possible to find linear combinations of both processes which are perfectly correlated. From a functional analysis standpoint, the covariance operators involved in the analysis require regularization as they are Hilbert-Schmidt and thus do not possess an inverse (see e.g., Rynne and Youngson \citep{Rynne01}).  The situation involved with functional canonical correlation analysis (FCCA) is analogous, therefore, to the classic inverse problem of finding approximate solutions to equations involving Freidholm integral equations.  Much like this classic problem, regularization plays an instrumental role in it's resolution.

This paper is organized as follows.  In Section 2 we introduce the notations, definitions and assumptions which we will utilize throughout the paper.  In this section we will also discuss why reproducing kernel Hilbert space (RKHS) methods is the ideal Hilbert space to solve the functional canonical correlation analysis (FCCA) problem. In Section 3 we will introduce the notions of canonical correlation and discuss why the Eubank and Hsing \citep{Eubank06} approach to FCCA provides most complete definition to canonical correlation analysis without regularization.  In Section 4 we will discuss the general theory associated with regularization and introduce both the Tikhinov and truncated singular value decomposition (TSVD) types of regularization (Engl et al. \cite{EHN00}).  In Section 5 and 6 we will discuss the consistency and asymptotic distributional theory associated with Tikhinov regularized canonical correlation operators, and in Sections 7 and 8 we will do the same with TSVD regularization.  Finally, Section 9 will be devoted to summarizing our conclusions and providing some further recommendations.

\section{Basic notation, definitions and assumptions}

Let $E$ be a subset of $\BBR$ and $\nu$ a sigma-finite measure on $E$. We then consider the case where a stochastic process $\{X(t), t \in E\}$ takes values in the Hilbert space $\CH = L^2(E)$ of square integrable functions on $E$ with inner product
$\la f,g \ra_{\CH} \equiv \int_{E} f(t)g(t) d \nu(t)$. Throughout it will be assumed that
\begin{equation}\label{4moment} \E \| X \|^{4}_{\CH} < \infty . \end{equation}
Under this assumption, $\E \la X, f \ra_{\CH} = \int_{\CH} \la x, f \ra_{\CH} dP(x) < \infty$ for all $f \in \CH$ with $P$ denoting the induced probability measure of $X$ on $\CH$.  The Riesz-Frechet representation theorem then ensures the existence of an element $\mu \in \CH$ such that $\E \la X, f \ra  = \la \mu,f \ra$.  Under assumption (\ref{4moment}), the Riesz-Frechet representation theorem also ensures the existence of the covariance operator $S: \CH \mapsto \CH$, which is given by
\begin{equation}\label{Sdef} \E[ \la f, X - \mu  \ra \la X - \mu , g \ra ] = \E[ \la f, (X - \mu) \otimes_{\CH} (X - \mu) g \ra ] = \la f, S g \ra \end{equation}
where $\otimes_{\CH}$ is the tensor product in $\CH$ and is defined by $( f \otimes_{\CH} g) h \equiv \la f, h \ra_{\CH} g$  for all $f, g, h \in \CH$.  We may also write $S = \E[ (X - \mu) \otimes_{\CH} (X - \mu) ]$.  It is also well known that the covariance operator $S$ is self-adjoint, non-negative definite, and has finite trace (see Laha and Rohatgi \citep{Rynne01}). The finite trace property ensures that $S$ is Hilbert-Schmidt and hence compact.

For any abstract Hilbert spaces $\CM$ and $\CN$ let $\CB(\CM, \CN)$ denote the Banach space of all bounded operators that map $\CM$ to $\CN$.  A subclass of $\CB(\CM, \CN)$ is $\CK(\CM, \CN)$ which will denote the set of all compact operators that map $\CM$ to $\CN$. Of particular importance in this paper is the subclass of compact operators which have finite trace, known as Hilbert-Schmidt operators. Let $\CK_{HS}(\CM, \CN)$  denote the set of all Hilbert-Schmidt operators that map $\CM$ to $\CN$. In this paper we will use the simplifying notation that $\CB(\CM) = \CB(\CM,\CM)$, $\CK(\CM) = \CK(\CM,\CM)$ and $\CK_{HS}(\CM) = \CK_{HS}(\CM,\CM).$ The ordinary operator norm on $\CB(\CM)$ will be denoted by $\| \cdot \|.$ The set of Hilbert-Schmidt operators $\CK_{HS}(\CM)$ becomes a separable Hilbert space when it is endowed with the inner product
\begin{equation}\label{HS inner product} \la A, B \ra_{HS} = \sum_{k=1}^{\infty} \la A e_{k}, B e_{k} \ra_{\CM} = \tr(A^{*}B), ~~ A,B \in \CK_{HS}(\CM)  \end{equation}
with $\{e_{k}\}_{k=1}^{\infty}$ denoting any complete orthonormal system (CONS) for $\CM$. This inner product does not depend on the
choice of basis (Kato \citep{Kato80}). The inner product, norm and tensor product on $\CK_{HS}(\CM)$ will be denoted by $\la \cdot, \cdot \ra_{HS}$, $\| \cdot \|_{HS}$ and $\otimes_{HS}$ respectively.

Next, assume that $\CH_1$ and $\CH_2$ are two closed subspaces of $\CH$ such that
\begin{equation*}  \CH = \CH_{1} \oplus \CH_{2},~~ \CH_1 \perp \CH_2 \end{equation*}
and let $\Upsilon_{i}$, $i = 1,2$ denote the orthogonal projection operator of $\CH$ onto $\CH_{i}$ for $i = 1,2.$ Suppose further that
$X_{i} = \Upsilon_{i} X$, $\mu_{i} = \Upsilon_{i} \mu$, and $S_{ij}$ denote the restriction of $S$ to $\CH_{i}$ and $\CH_{j}$ for $i,j = 1,2$ so that   $S_{ij} = \Upsilon_{j} S \Upsilon_{i}$.  Because the $\Upsilon_{i}$ are bounded and $S$ is Hilbert-Schmidt, the $S_{ij}$ for $i,j=1,2$ are also Hilbert-Schmidt and compact. In addition, the $S_{ii}$ are self-adjoint and non-negative definite. For convenience, we henceforth denote $S_{ii} = S_{i}.$

For $i = 1,2$, let $\{\phi_{in}\}_{n=1}^{\infty}$ be an orthonormal basis corresponding to eigenvectors of $S_{i}$ with
$\{\lambda_{in}\}_{n=1}^{\infty},$ the corresponding sequence of non-negative
eigenvalues. Since $S_{i}$ is self-adjoint, non-negative and compact we may write
\begin{equation}\label{S otimes expansion} S_{i} = \sum_{n=1}^{\infty} \lambda_{in} \phi_{in} \otimes_{\CH_i} \phi_{in}, ~~ i = 1,2 \end{equation}
with $\lambda_{i1} \geq \lambda_{i2} \geq \cdots \geq 0$ a decreasing sequence whose only limit can be zero. For our purposes we might as well assume without loss of generality (WLOG) that $\{\phi_{in}\}_{n=1}^{\infty}$ is a CONS for $\CH_{i}$, $S_{i}$ is strictly positive and $\CH_{i} = \ker(S_{i})^{\perp}$ for $i=1,2.$ We make this assumption since if $\varphi \in \ker(S_{i})$ then $\var[ \la \varphi,X_{i} \ra_{\CH}] = \la \varphi,S_{i} \varphi \ra_{\CH} = 0$, which would have the consequence that $\la \varphi,X_{i} \ra_{\CH_{i}} = \la \mu_{i}, \varphi \ra_{\CH_{i}}$ with probability one. It is also convenient at this juncture to assume, WLOG, that the mean of the process is zero because if this does not hold we may always consider the covariance of the process $X(\cdot) - \mu(\cdot)$ instead. It should be mentioned that in (\ref{S otimes expansion}) the list of eigenvalues $\{\lambda_{in}\}$ is repeated according to their multiplicity. An alternative expression for (\ref{S otimes expansion}) involving eigenprojection operators is
\begin{equation}\label{T expressed as a projection operator} S_{i} = \sum_{h=1}^{\infty} \tilde{\lambda}_{ih} P_{ih}, ~~ \text{ for } i = 1,2 \end{equation}
where $\{\tilde{\lambda}_{ih}\}$ are the distinct elements of $\{\lambda_{in}\}$, and $P_{ih}$ is the finite dimensional projection operator onto the eigenspace associated with each distinct $\tilde{\lambda}_{ih}$ given by
\begin{equation}\label{Form for the orthogonal projection operator} P_{ih} = \sum_{\phi_{in}: ~ \lambda_{in} = \tilde{\lambda}_{ih}} \phi_{in} \otimes_{\CH_i} \phi_{in}. \end{equation}

Since the processes $\{X_{i}(\cdot)\}_{i=1}^{2}$ are of second order, they admit a Karhunen\textendash Lo\`{e}ve expansion
$X_{i}(\cdot) = \sum_{n=1}^{\infty} Z_{in}  \phi_{in}(\cdot)$ with the random variables $Z_{in}$ defined by $Z_{in} = \la X_{i},\phi_{in} \ra_{\CH_{i}}$ (see Ash and Gardiner \citep{Ash75} or Doob \citep{Doob53}). These variables are orthogonal in the sense that $\cov[ Z_{ij}, Z_{ik}] = \la \phi_{ij}, S_i \phi_{ik} \ra_{\CH_i} = \lambda_{ij} \delta_{jk}$ with $\delta_{jk}$ denoting the Kronecker delta function. Mercer's theorem then ensures that the covariance functions of the processes $\{X_{i}(t), t \in E_i\}_{i=1}^{2}$ are
\begin{align}\label{covariance kernel} K_{ii}(s,t) &= \E[X_{i}(s)X_{i}(t)] = \sum_{n=1}^{\infty} \sum_{m=1}^{\infty} \E[Z_{im} Z_{in}] \phi_{in}(s) \phi_{im}(t) \notag\\ &= \sum_{n=1}^{\infty} \lambda_{in} \phi_{in}(s) \phi_{in}(t) ~ \text{ for } i=1,2. \end{align}
Moreover, the cross-covariance kernel is then
\begin{align}\label{cross-covariance kernel} K_{12}(s,t) &= \E[ X_{1}(s) X_{2}(t)]  = \sum_{n=1}^{\infty} \sum_{m=1}^{\infty} \E[Z_{1n} Z_{2m}] \phi_{1n}(s) \phi_{2m}(t) \notag\\ &= \sum_{n=1}^{\infty} \sum_{m=1}^{\infty} \gamma_{mn} \phi_{1n}(s) \phi_{2m}(t) \end{align}
and we note that $K_{12}(s,t) = K_{21}(s,t)$.  For notational simplicity let $K_{ii}(s,t) = K_{i}(s,t)$. It is also well known that for all $f \in \CH_{i}$
\begin{equation}\label{S is a integral operator} (S_{ij} f)(t) = \int_{E} K_{ij}(s,t) f(s) d\nu(s) = \la K_{ij}(\cdot,t),f \ra_{\CH_i} ~ \text{ for } i, j = 1,2. \end{equation}
An alternate form for the cross-covariance operator $S_{12}: \CH_2 \mapsto \CH_1$ is given by
\begin{equation} S_{12} = \sum_{n=1}^{\infty} \sum_{m=1}^{\infty} \gamma_{mn} \phi_{2n}(s) \otimes_{\CH_2} \phi_{1m}(t) = S^{*}_{21}. \end{equation}

In addition to $\CH = \ker(S)^{\perp} \subseteq L^2(E)$, two additional types of Hilbert spaces will play prominent roles in further developments. The first type of Hilbert space are the reproducing kernel Hilbert spaces (RKHS) associated with the symmetric covariance kernels $K_{i}(\cdot,*)$, denoted $\CH(K_i)$ (see Aronszajn \citep{Aronszajn50} or Berlinet and Thomas-Agnan \citep{Berlinet04}).  The second type are the Hilbert spaces generated by each stochastic process, denoted $L^2_{X_i}$ with $i=1,2$ (see Parzen \citep{Parzen61}). To construct both of these Hilbert spaces we first let $\{t_1,\ldots,t_n\}$ be any finite collection of points in $E$ and let $\BX_{in} = \left[ X_{i}(t_{1}),\ldots, X_{i}(t_n) \right]'$ with $\BK_{in} = \{K_{i}(t_j,t_k)\}_{j,k=1}^n$ denoting the covariance matrix of $\BX_{in}$ for each $n \in \BBN$. Next we define the pre-Hilbert space generated by the process to be the set of all arbitrary finite dimensional linear combinations of the process, i.e. $L^2_{X_{in}} = \{U= \ba' \BX_{in} \text{ : } \ba \in \ker(\BK_{in})^{\perp} \subseteq \BBR^n \}$ where the inner product between two elements is given by
\begin{equation}\label{inner product on L2X} \la \ba'\BX_{in},\bb'\BX_{in} \ra_{L^2_{X_{i}}} = \cov[\ba'\BX_{in},\bb'\BX_{in}] = \ba' \BK_{in} \bb.\end{equation}
Likewise, the pre-Hilbert space of the RKHS is defined to be the column space of $\BK_{in}$, i.e. $\CH(\BK_{in}) = \{ \bbf = \BK_{in} \ba \text{ : } \ba \in \ker(\BK_{in})^{\perp} \subseteq \BBR^n \}$ and the inner product between any $\bbf = \BK_{in} \ba$ and $\bg = \BK_{in} \bb$ is given by
\begin{equation}\label{inner product on CHK} \la \boldf, \boldg \ra_{\CH(\BK_{in})} \equiv \bbf' \BK_{in}^{\dag} \bg = \ba'\BK_{in}\bb \end{equation}
where $\BK_{in}^{\dag}$ denotes the Moore-Penrose inverse of $\BK_{in}$.  The Parzen--Lo\'{e}ve congruence mapping is determined uniquely by $\Psi_{i}(K_{i}(t,\cdot)) = X_{i}(t)$ for each $t \in E$ with the result that every linear combination $U$ of the $\BX$ vector with nonzero variance can be expressed as
\begin{equation}\label{Parzen congruence} U = \Psi(\bbf) =  \bbf' \BK_{in} \BX_{in} \end{equation}
for some $\bbf \in \CH(\BK_{in})$ (see King \citep{King10}).  It is a simple matter to see that inner product given by (\ref{inner product on CHK}) satisfies the reproducing property (see Aronszajn \citep{Aronszajn50}). To see this let $k$ be any index in $1,\ldots, n$ and let $\BK_{in}(t_k,\cdot) = \BK^{T}_{in}(\cdot, t_k)$ denote the $k^{th}$ row of $\BK_n$. Now for any $\bbf = \BK_{in} \ba \in \CH(\BK_{in})$  we have that
$$ \la \BK_n(\cdot,t_k),\boldf \ra_{\CH(\BK_{in})} = \BK_{in}(t_{k},\cdot) \boldK_{in}^{\dag} \boldK_{in} \bolda = \boldK_{in}(t_{k},\cdot) \bolda = \bbf(t_k).$$
This demonstrates that the pre-Hilbert space $\CH(\boldK_n)$ given by inner product defined in (\ref{inner product on CHK}) must be the unique RKHS of the process $\{X(t_i)\}_{i=1}^{n}$. To complete the construction of $\CH(K_i)$ and $L^2_{X_i}$ we then extend the realm of the pre-Hilbert spaces, which presently apply to any finite collection of points $\{t_1,\ldots,t_n\} \in E$ to the index set, $E$ in its entirety. This construction is accomplished through Cauchy completion or adding in the limits of arbitrary linear combinations of the form $\sum_{j=1}^{n} a_i K(\cdot,t_{j})$ and  $\sum_{j=1}^{n} a_i X(t_j)$.  In this fashion, we see that
\begin{equation}\label{HK prehilbert closure} \CH(K_{i}) = \overline{\{f: f(\cdot) = \int_{E} K(\cdot,t) a(t) d\nu(t)\} } = \overline{ \Ima(S_{i})} = \ker(S_{i})^{\perp}\end{equation}
and
\begin{equation}\label{L2X prehilbert closure} L^2_{X_i} = \overline{\{U: U = \int_{E} a(t) X(t) d\nu(t)\} } \end{equation}
with $\overline{A}$, denoting the closure of any set $A$. In this infinite dimensional setting the RKHS is the set of function on $E$ given by
\begin{equation}\label{RKHS Representation} \CH(K_{i}) = \{f: f(\cdot) = \sum_{j=1}^{\infty} \lambda_{ij} f_{ij} \phi_{ij}(\cdot), \|f\|^2_{\CH(K_{i})} = \sum_{j=1}^{\infty} \lambda_{ij} f_{ij}^2 < \infty\} \end{equation}
where $f_{ij} = \la f, \phi_{ij} \ra_{\CH_{ij}}$ are the generalized Fourier coefficients relative to the CONS $\{\phi_{ij}\}_{j=1}^{\infty}$ for $\CH_{i}$, $i=1,2$.  An application of the integral representation theorem of Parzen \cite{Parzen61} then produces the following result.

\begin{theorem}\label{Congruence Map Psi to L2X}
For $i = 1,2$ let $f(\cdot) = \sum_{j=1}^{\infty} \lambda_{ij} f_{ij} \phi_{ij}(\cdot)$ be in $\CH(K_{i})$. Then,
\begin{equation}\label{Psi of f} \Psi_{i}(f) = \sum_{j=1}^{\infty} f_{ij} Z_{ij}   \text{ and } \Psi_{i}^{-1}\left( \sum_{j=1}^{\infty} f_{ij} Z_{ij} \right) = \sum_{j=1}^{\infty} \lambda_{ij} f_{ij} \phi_{ij} \end{equation}
with $Z_{ij} = \la X_{i}, \phi_{ij} \ra_{\CH_{i}}$ and $\Psi_{i}^{-1} = \Psi_{i}^{*}$, where $\Psi_{i}^{*}$ denotes the adjoint of $\Psi_{i}$.
\end{theorem}

The importance of the RKHS inner product when formulating theory regarding integral operators was shown by Nasheed and Wahba \citep{Wahba74}.  These authors provided a characterization of the RKHS $\CH(K_{i})$ generated by the kernel $K_{i}$ to closure of the image of the integral operator for the symmetric square root, $\overline{\Ima(S_{i}^{1/2})}$. In this regard, first notice that since the $S_{i}$ are positive (and self-adjoint), they have symmetric square roots $S_{i}^{1/2}$ with associated symmetric kernel $\Phi_{i}(s,t)$ given explicitly by
\begin{equation}\label{Defn of kernel Phi} \Phi_{i}(s,t) = \sum_{j=1}^{\infty} \lambda_{ij}^{1/2} \phi_{ij}(s) \phi_{ij}(t), ~~ i = 1,2.\end{equation}
For $i=1,2$ the symmetric kernels $\Phi_{i}(s,t)$ satisfy
\begin{equation}\label{K in terms of Q the kernel} K_{i}(s,t) = \int_E \Phi_{i}(s,r)\Phi_{i}(t,r) d\nu(r).\end{equation}
We further note that
\begin{equation*}\label{ker S perp relation} \Ima(S_{i}^{1/2}) \subseteq \overline{\Ima(S_{i}^{1/2})} = \ker(S_{i}^{1/2})^{\perp} =  \ker(S_{i})^{\perp} = \CH_{i}. \end{equation*}
Nasheed and Wahba \citep{Wahba74} then arrive at the following important theorem.

\begin{theorem}(Nasheed and Wahba, 1974)\label{Nasheed Wahba 1}
For $i=1,2$ the RKHS $\CH(K_{i})$ consist of functions of the form
\begin{equation*}\label{f as innprod with Q} f(\cdot) = \int_E g(s) \Phi_{i}(\cdot,s) d\nu(s)\end{equation*}
for some $g \in \CH_{i}$.  The inner product in $\CH(K_{i})$ is
\begin{equation}\label{congruence2} \la f_1,f_2 \ra_{\CH(K_{i})} = \la g_1,g_2 \ra_{\CH_{i}}\end{equation}
where $g_1, g_2 \in \CH_{i}$ are the minimal $L^2(E)$ norm solutions of
\begin{equation*}\label{f as innprod with Q2} f_j(\cdot) = \int_E g_j(s) \Phi_{i}(\cdot,s) d\nu(s), ~ j = 1,2.\end{equation*}
\end{theorem}

\noindent Proof: For $i = 1,2$, let $V_{i}$ be the smallest closed subspace of $\CH_{i}$
that contains $\Phi_{i}(\cdot,t)$ for all $t \in E$. Since the smallest linear space containing $\Phi_{i}(\cdot,t)$ for all $t \in E$
is ${\rm span} \{\Phi_{i}(\cdot,t) : t \in E\}$, it follows that $V = \overline{{\rm span}} \{\Phi_{i}(\cdot,t) : t \in E\}$. Now the projection theorem ensures that for each $f \in \CH(K)$ there exists a unique element $g_f \in V$ of minimal norm which is the best-approximate solution to the inverse problem
\beqs\label{interpolation} f(t) = (S_{i}^{1/2} g_f)(t) = \int_E g_f(s) \Phi_{i}(s,t) d\nu(s), ~ \forall t \in E.\eeqs
Because $g_f$ is unique, the inner product given by (\ref{congruence2}) and associated
norm are well defined.  We now only need to show that $K_{i}(\cdot,\cdot)$ are the reproducing kernel. However,
\beqs\label{K in terms of Q2} K_{i}(\cdot,t) = \la \Phi_{i}(*,\cdot),\Phi{i}(*,t) \ra_{\CH_{i}}.\eeqs
Thus, by (\ref{congruence2}),
\beqs\label{reproducing proof} \la K_{i}(\cdot,t),f \ra_{\CH(K_{i})} = \la \Phi_{i}(*,t),g_f(*) \ra_{\CH_{i}} = (S_{i}^{1/2} g_f)(t) = f(t) .\eeqs \qed

This theorem shows that for $i =1,2$ the optimal Hilbert space to solve inverse problems associated with integral equations of the form $f = S_i^{1/2} g$ is in the RKHS setting $\CH(K_{i})$. To illustrate this, consider the problem of finding a function $g(\cdot)$ to satisfy $S_i^{1/2}g = \int_E g(s) \Phi_{i}(\cdot,s) d\nu(s) = f(\cdot)$ for some given $f(\cdot) \in \CH_i$.  A least-squares solution to this problem is a minimizer of $\|S_i^{1/2} g - f\|_{\CH_i}$ and a best least-squares solution is the one with minimum norm. If we let $F_i = (\Ima S_i^{1/2}) \oplus (\Ima S_i^{1/2})^{\perp}$ and assume $f \in F_i$, $g$ is a least squares solution if and only if $ S_i^{1/2} S_i^{1/2} g = S_i^{1/2} f = S g.$
Furthermore, the unique best least-squares solution is given by $g = (S_i)^{\dag} S_i^{1/2} f = S_i^{1/2 \dag} f$
with $S_i^{1/2 \dag}$ denoting the Moore-Penrose inverse of $S_i^{1/2},$ and no least-squares solution exists if $f \notin F_i$.
However, from Engl et al. \citep{EHN00}, $f \in F_i$ if and only if $f$ satisfies the Picard criterion
\begin{equation}\label{picardQ} \sum_{j=1}^{\infty} \frac{\la f,\phi_{ij} \ra_{\CH_i}^2}{\lambda_{ij}} < \infty \end{equation}
and, in that case,
\beqs\label{picards condition characterization} g(\cdot) = (S_{i}^{1/2 \dag} f)(\cdot) =  \sum_{j=1}^{\infty} \frac{ \la f,\phi_{ij} \ra_{\CH_i}}{\sqrt{\lambda_{ij}}} \phi_{ij}(\cdot).\eeqs
Note that $f \in F_i$ if and only if $\| f \|_{\CH(K_i)} = \sum_{j=1}^{\infty} \frac{\la f,\phi_{ij} \ra_{\CH_i}^2}{\lambda_{ij}} < \infty.$ Consequently, $\CH(K_i) = \overline{ (\Ima S_{i}^{1/2} )} = \ker(S_i)^{\perp},$ under the inner product
\beqs\label{last characterization} \la f_{i1},f_{i2} \ra_{\CH(K_i)} = \la S_{i}^{1/2 \dag} f_{i1}, S_{i}^{1/2 \dag} f_{i2} \ra_{\CH_i} \eeqs
with
\beqs\label{Qdag on f} S_{i}^{1/2 \dag} f_{ij}(\cdot) = \sum_{k=1}^{\infty} \frac{ \la f_{ij},\phi_{ik} \ra_{\CH_i}}{\sqrt{\lambda_{ik}}} \phi_i(\cdot)\eeqs
for $i,j = 1,2$.

For further developments, a congruence which connects $\CH_{i}$ to $\CH(K_{i})$ must be established.

\begin{corollary}(Eubank and Hsing, 2008)\label{congruence thm1}
For $i=1,2$ the Hilbert spaces $\overline{\Ima (S_{i}^{1/2})} = \ker(S_{i})^{\perp}$ and $\CH(K_{i})$ are congruent
under the mapping $\Gamma_{i}: \CH_{i} \mapsto \CH(K_{i})$ defined by
\begin{equation} \label{gamma} (\Gamma_{i} g)(\cdot) \equiv \sum_{j=1}^{\infty} \sqrt{\lambda_{ij}}g_{ij} \phi_{ij}(\cdot) \end{equation}
where $g = \sum_{j=1}^{\infty} \la g, \phi_{ij} \ra_{\CH_{i}} \phi_{ij} = \sum_{j=1}^{\infty} g_{ij} \phi_{ij} \in \ker(S)^{\perp}$.
The inverse mapping
\begin{equation}\label{gamma inv} (\Gamma_{i}^{-1} f)(\cdot) \equiv \sum_{j=1}^{\infty} \sqrt{\lambda_{ij}} f_{ij} \phi_{ij}(\cdot) \end{equation}
for $f = \sum_{j=1}^{\infty} \lambda_{ij} f_{ij} \phi_{ij}(\cdot) \in \CH(K_{i}),$ is also the adjoint of $\Gamma_{i}$.
\end{corollary}

Note that for $i = 1,2$ the operators $\Gamma_i$ and $S^{1/2}_{i}$ are equal in the sense that for any $f \in \CH_i$, $\Gamma f = \sum_{j=1}^{\infty} \sqrt{\lambda_{ij}} f_{ij} \phi_{ij} = S^{1/2} f$.  The difference is in terms of the norm and inner product for
the range of each operator.

\section{Canonical Correlation}

The literature on functional canonical correlation can be roughly dichotomized
into formulations involving Hilbert space valued processes in $\CH_{i} = \ker(S_i)^{\perp}$ (see He et al. \citep{He00} \citep{He02}) and an
alternative approach that relies on reproducing kernel Hilbert space (RKHS) theory (Eubank and Hsing, \citep{Eubank06}). In this section we will compare and contrast these two different approaches to functional CCA. In the He et al. \citep{He00} approach the $k^{th}$ squared
canonical correlation $\rho^2_{k}$ and associated weight functions $f_{k}$ and $g_{k}$ are found by the singular value decomposition
of the cross-correlation operator of $X_1$ and $X_2$ defined by
\begin{equation}\label{Cross correlation operator} R = S_1^{1/2 \dag} S_{12} S_2^{1/2 \dag} \end{equation}
where $S_{i}^{1/2 \dag}$ denotes the Moore-Penrose generalized inverse of $S_{i}^{1/2}$ for $i = 1,2$ and
is given explicitly by
\begin{equation}\label{ginv of S} S_{i}^{1/2 \dag} = \sum_{h=1}^{\infty} \tilde{\lambda}^{-1/2}_{ih} P_{ih}. \end{equation}
The right hand and left hand eigenvectors are then found by eigenvalue and eigenvector analysis of the operators $RR^{*}$ and $R^{*}R$. The basic problem is that unlike the usual situation in the finite-dimensional case, the square roots of covariance operators of infinite dimensional Hilbert space valued processes are not invertible. To resolve this issue, He et al. \citep{He00} restricts the domain of $\{\CH_1, \CH_2\}$ to the subspace where the Moore-Penrose inverses of $S_1^{1/2}$ and $S_2^{1/2}$ can be defined. Thus, for $i = 1,2$, the domain of $S_i^{1/2 \dag}$ is restricted to $F_{i} \equiv \{S_i^{1/2} h : h \in \ker(S_i)^{\perp} \}$ and is characterized as the set of functions satisfying the Picard criterion (\ref{picardQ}) (see Engl et al. \citep{EHN00}). Now, subject to the restriction that the domain of $R$ be $F_2$, let $\rho^2_1 \geq \rho^2_2 \geq \cdots \geq 0$ denote the eigenvalues of $R^{*}R$ with $g_{1}, g_{2},\ldots \in F_2$ the corresponding eigenvectors. The left hand eigenvectors are obtained by $f_{k} = R g_{k}/ \rho_k \in F_1$. The canonical correlations and weight functions are $\{ \rho_{k}, u_{k} = S^{1/2 \dag}_{1} f_{k}, v_{k} = S^{1/2 \dag}_{2} g_{k}\}_{k=1}^{\infty}$ and the corresponding canonical variables are $\{ U_{k} = \la u_{k}, X_1 \ra_{\CH_1}, V_{k} = \la v_k, X_2 \ra_{\CH_2} \}_{k=1}^{\infty}$. In the He et al. \citep{He00} method the weight functions are not well defined whenever $f_k \notin F_1$ or $g_k \notin F_2$.

In contrast to the He et al. \citep{He00} method, the approach of Eubank and Hsing \citep{Eubank06} involves the singular value decomposition
of the RKHS based operator $T:  \CH(K_2) \mapsto \CH(K_1)$ defined such that for any $\tilde{g} \in \CH(K_2)$
\begin{equation}\label{T defn} (T \tilde{g})(s) = \la K_{12}(s,\cdot), \tilde{g} \ra_{\CH(K_{2})}. \end{equation}
Let $\{ \rho^2_{k}, \tilde{g}_{k} \}_{k=1}^{\infty}$ denote the eigenvalues and eigenvectors of $T^{*}T$ and $\{ \rho^2_{k}, \tilde{f}_{k} \}_{k=1}^{\infty}$ denote the eigenvalues and eigenvectors for $TT^{*}.$ Then
\begin{equation}\label{T decomposition} T = \sum_{k=1}^{\infty} \rho_{i} \tilde{g}_k \otimes_{\CH(K_2)} \tilde{f}_k. \end{equation}
The $k^{th}$ canonical correlation is $\rho_k$ and $\{ \tilde{f}_k, \tilde{g}_k\}_{k=1}^{\infty}$ are
the canonical weight vectors in $\{ \CH(K_1), \CH(K_2) \}$. These canonical weight vectors correspond to
the canonical variables $\{\Psi_{1}(\tilde{f}_j), \Psi_{Y}(\tilde{g}_j) \}_{j=1}^{\infty}$ that
represent the maximally correlated elements of $\{L^2_{X_1}, L^2_{X_2}\}$.

The relationship between $T$ and $R$ was established by Eubank and Hsing \citep{Eubank06} and can be simply derived by substituting the expression for $K_{12}$ given by (\ref{cross-covariance kernel}) into (\ref{T defn}). It follows that for any $\tilde{g} \in \CH(K_2)$,
\begin{align}\label{T expansion right after} (T \tilde{g})(s) &= \sum_{j=1}^{\infty} \sum_{k=1}^{\infty} \gamma_{jk}  \la \phi_{2k} , \tilde{g} \ra_{\CH(K_{2})} \phi_{1j}(s) \notag\\ &= \sum_{j=1}^{\infty} \sum_{k=1}^{\infty} (\rho_{jk} \sqrt{\lambda_{1j} \lambda_{2k}}) \la \phi_{2k} , \tilde{g} \ra_{\CH(K_{2})} \phi_{1j}(s) \end{align}
with $\rho_{jk} = \frac{\gamma_{jk}}{\sqrt{\lambda_{1j} \lambda_{2k}}}$. Now, since $\{\phi_{ik}\}_{k=1}^{\infty}$ are CONSs for $\ker(S_i)^{\perp}$, it follows that $\{\tilde{\phi}_{ik} = \Gamma_i \phi_{ik} = \sqrt{\lambda_{ik}}\phi_{ik}\}_{k=1}^{\infty}$ are CONSs for $\CH(K_i)$.
As a result, the operator $T$ may be written as
\begin{eqnarray}\label{T3} T &=& \sum_{j=1}^{\infty} \sum_{k=1}^{\infty} \rho_{jk} \tilde{\phi}_{2k} \otimes_{\CH(K_2)} \tilde{\phi}_{1j} \nonumber\\ &=& \sum_{j=1}^{\infty} \sum_{k=1}^{\infty} \rho_{jk} [(\Gamma_2 \phi_{2k} ) \otimes_{\CH(K_2)} (\Gamma_1 \phi_{1j})]
\nonumber\\ &=& \sum_{j=1}^{\infty} \sum_{k=1}^{\infty} \rho_{jk} [(\phi_{2k} \Gamma^{-1}_2) \otimes_{\CH_2} (\Gamma_1 \phi_{1j})]
\nonumber\\ &=&   \Gamma_1 \left[ \sum_{j=1}^{\infty} \sum_{k=1}^{\infty} \rho_{jk} [\phi_{2k} \otimes_{\CH_2} \phi_{1j}] \right] \Gamma^{-1}_2
\nonumber\\ &=&   \Gamma_1 R \Gamma^{-1}_2 = \Gamma_1 S_1^{1/2 \dag} S_{12} S_2^{1/2 \dag} \Gamma^{-1}_2. \end{eqnarray}
Because $\Gamma_2$ is a bijection, $\Gamma^{-1}_1 T \Gamma_2: \CH_2 \mapsto \CH_1$ has the form
\beqs\label{T4} \Gamma^{-1}_1 T \Gamma_2 = \sum_{j=1}^{\infty} \sum_{k=1}^{\infty} \rho_{jk} \left( \phi_{2k} \otimes_{\CH_2}  \phi_{1j} \right) \eeqs and the domain is $\ker(S_2)^{\perp}$. By contrast, if we utilize the He et al. \citep{He00} method and restrict ourselves to the domain $F_{2} = \Ima(S_2^{1/2}) \subset \overline{\Ima(S_2^{1/2})} = \ker(S_2)^{\perp}$ then, on this restricted subspace of $\ker(S_2)^{\perp}$,
\begin{equation*}\label{matrix rep of R} R = S_1^{1/2 \dag} S_{12} S_2^{1/2 \dag} = \sum_{j=1}^{\infty} \sum_{k=1}^{\infty} \rho_{jk} \left( \phi_{2k} \otimes_{\CH_2} \phi_{1j} \right) = \Gamma^{-1}_1 T \Gamma_2 |_{F_2}. \end{equation*}
Since $\Gamma_1$ and $\Gamma_2$ are unitary, $T$ is unitarily equivalent to $R$ and the two methods agree when both methods are well defined.  The differences between the approaches can be briefly summarized by the fact that in the He et al. \citep{He00} approach the domain of $R$ must be restricted to $F_2 = S_2^{1/2} \left(\ker(S_2)^{\perp}\right)$, which is a dense proper subset of $\ker(S_2)^{\perp}$ in the infinite dimensional case. By contrast, the domain of $\Gamma^{-1}_1 T \Gamma_2$ in Eubank and Hsing \citep{Eubank06} approach is all of $\ker(S_2)^{\perp}$, since the mapping $\Gamma_2$ is a unitary bijective mapping from $\ker(S_2)^{\perp} \mapsto \CH(K_2)$. Therefore the Eubank and Hsing \citep{Eubank06} approach is the more comprehensive definition while the He et al. \citep{He00} approach can have non-attainable solutions on the boundary $\left( \overline{\Ima(S_2^{1/2})} \setminus \Ima(S_2^{1/2}) \right)$.  This reveals the advantage of RKHS based formulation and we will therefore consider asymptotics associated with the regularized approximations to $TT^{*}$ and $T^{*}T$ rather than $RR^{*}$ and $R^{*}R$ in this paper.

\section{Regularization}

The need to employ some form of regularization in the functional data analysis setting is well established on both theoretical as well as computational grounds by many authors.  For example, it was perhaps Leurgans et al. \citep{lms93} who first observed that the sample
covariance operator of a stochastic process has a finite dimensional kernel (Riesz \& Sz.-Nagy \citep{Riesz78}), while acting on an infinite dimensional space. Cupidon et al. \citep{Cupidon06} then showed how most of the deficiencies of the population canonical correlation can be remedied if a regularized approximation to the inverses of the covariance operators are involved.

If $B \in \CK(\CH_1,\CH_2)$ is arbitrary and we are given $g \in \CH_2$, it often happens that we are asked to solve the equation $B f = g$.
If $\{\beta_n, \phi_n, \theta_n\}_{n=1}^{\rank(B^{*}B)}$ is the singular system for $B$ so that
\beqs\label{form of B} B = \sum_{n=1}^{\rank(B^{*}B)} \beta_n \phi_n \otimes_{\CH_1} \theta_n \eeqs
and $g \in \Ima(B) \oplus \Ima(B)^{\perp}$, then it is well known that a unique best approximate (least squares)
solution $f_{*}$ exists and is given by
\beqs\label{moore penrose solution of T} f_{*} =  B^{\dag} g =  \sum_{n=1}^{\rank(B^{*}B)} \frac{\la y,\theta_n \ra_{\CH_2}}{\beta_n} \phi_n.\eeqs
For a compact operator $B$, $B f = g$ is often ill-posed (see e.g., Theorem 2.14 of Vogel \cite{Vogel02}) and attempts to directly
use $B^{\dag}$ will result in numerically unstable algorithms.  The standard approach to dealing with this problem
is to replace $B^{\dag}$ with a family of so called regularization operators $D(\alpha): \CH_2 \mapsto \CH_1$ that are
indexed by a regularization parameter, $\alpha \in (0,a) \subset \BBR$, with $a > 0$.  The family $\{ D(\alpha) : \alpha \in (0,a)\}$
approximates $B^{\dag}$ in the sense of the following definition (see Vogel \cite{Vogel02} p. 22-23).

\begin{definition}\label{definiton of regularization}
The family $\{D(\alpha) : \alpha \in (0,a) \}$ is a regularization scheme which converges to $B^{\dag}$ if
\begin{enumerate}
\item for each $\alpha \in (0,a)$, $D(\alpha)$ is a continuous operator and
\item given any $g \in \Ima(B)$, for any sequence $\{g_n\} \subset \CH_2$ which converges to $g$,
one can pick a sequence $\{\alpha_n\} \subset (0, a)$ such that
\beqs\label{formula for the regularized solution} \left[D(\alpha_n)\right](g_n) \rightarrow  B^{\dag} g ~~ \text{ as } n \rightarrow \infty. \eeqs
\end{enumerate}
\end{definition}

\noindent Of particular interest are linear regularization schemes which have singular value representations as
\beqs\label{singular value rep of R} D(\alpha) =  \sum_{n=1}^{\rank(B^{*}B)} \frac{w_{\alpha}(\beta_n^2)}{\beta_n} \theta_n \otimes_{\CH_2} \phi_n \eeqs
where $w_{\alpha}(\beta_n^2)$ is a real valued function of the squared singular values and
$\alpha$ is such that $w_{\alpha}(\beta_n^2) \rightarrow 1$ as $\alpha \rightarrow 0$. The function
$w_{\alpha}(\beta_n^2)$ is called the filter function (see Engl et al. \citep{EHN00}). Two of the most popular examples for filters are
\begin{equation}\label{Tikhinov Filter} w_{\alpha}(\beta_n^2) =  \frac{ |\beta_n|}{ |\beta_n| + \alpha},
~~ \text{ for } \alpha \in (0,\infty) \text{ and } n = 1, \ldots, \rank(B^{*}B)  \end{equation}
and
\begin{equation}\label{TSVD Filter} w_{\alpha}(\beta_n^2) =  \left\{ \begin{array}{cc} 1 & \text{ if } \beta_n^2 > \alpha
\\ 0 & \text{ if } \beta_n^2 \leq \alpha \end{array} \right.
 \text{ for } \alpha \in (0,\|B\| ] \text{ and } n = 1, \ldots, \rank(B^{*}B). \end{equation}
Equation (\ref{Tikhinov Filter}) is referred to as the Tikhinov filter function and (\ref{TSVD Filter}) is referred to as
the truncated singular value decomposition (TSVD) filter function.  In the case of TSVD regularization,
the parameter $\alpha$ in (\ref{TSVD Filter}) determines the cut-off or threshold level for the TSVD regularization and produces
\begin{equation*}\label{singular value rep of TSVD reg op} D(\alpha) =  \sum_{\beta_n^2 > \alpha} \beta_n^{-1} \theta_n \otimes_{\CH_2} \phi_n \end{equation*}
which is a finite rank operator whenever $\alpha > 0$. This paper will focus on asymptotics associated with Tikhinov and TSVD regularization schemes.
In developments which follow a Hilbert-Schmidt operator $B: \CH_1 \mapsto \CH_2$ will often be expressed in the form
\begin{equation}\label{arbirary form of hs operator} B =  \sum_{j=1}^{\infty} \sum_{k=1}^{\infty}  c_{jk} \phi_j \otimes_{\CH_1} \theta_k \end{equation}
with $\{\phi_j\}$, $\{\theta_k\}$ CONSs for $\CH_1$ and $\CH_2$, respectively. As $B \in \CK_{HS}(\CH_1, \CH_2)$ the
coefficients $c_{jk} \in \BBR$ will satisfy
\begin{equation*}\label{Hilbert schmidt norm of B} \|B\|^2_{HS} =  \sum_{j=1}^{\infty} \la B \phi_j, B \phi_j \ra_{\CH_1}
= \sum_{j=1}^{\infty} \| \sum_{k=1}^{\infty} c^2_{jk} \theta_k\|_{\CH_2} =  \sum_{j=1}^{\infty}\sum_{k=1}^{\infty} c^2_{jk} < \infty. \end{equation*}
The regularized operator $B(\alpha): \CH_1 \mapsto \CH_2$ will also be a Hilbert-Schmidt operator and have the form
\begin{equation}\label{a form for regularized hs operator} B(\alpha) =  \sum_{j=1}^{\infty} \sum_{k=1}^{\infty}  c_{jk}(\alpha) \phi_j \otimes_{\CH_1} \theta_k \end{equation}
with $\|B(\alpha)\|^2_{HS} = \sum_{i,j=1}^{\infty} c^2_{ij}(\alpha) < \infty$ and the property
$c_{jk}(\alpha) \rightarrow c_{jk}$ as $\alpha \downarrow 0$. As a result, the following theorem holds.

\begin{theorem}\label{HS Convergence theorem}
Suppose $B, B(\alpha) \in  \CK_{HS}(\CH_1, \CH_2)$ are of the forms (\ref{arbirary form of hs operator}) and (\ref{a form for regularized hs operator}),
respectively. If $c_{jk}(\alpha) \rightarrow c_{jk}$ as $\alpha \rightarrow 0$, then $\|B(\alpha) - B\| \rightarrow 0$.
\end{theorem}

\noindent Proof:  Let $A(\alpha) = B(\alpha) - B$ and $a_{jk}(\alpha) = (c_{jk}(\alpha) - c_{jk})$ so that
\beqs\label{a form for A alpha} A(\alpha) =  \sum_{j=1}^{\infty} \sum_{k=1}^{\infty}  a_{jk}(\alpha) \phi_j \otimes_{\CH_1} \theta_k. \eeqs
Since $A(\alpha)$ is Hilbert-Schmidt,
$\|A(\alpha)\|^2_{HS} = \sum_{j,k=1}^{\infty} a^2_{jk}(\alpha) < \infty$ for all permissable values of
the regularization parameter $\alpha$.  Consequently,
\begin{align}\label{Limit of HSnorm of A} \lim_{\alpha \rightarrow 0} \|A(\alpha)\|^2_{HS} &= \lim_{\alpha \rightarrow 0} \sum_{j,k=1}^{\infty}  a^2_{jk}(\alpha)
\notag\\ &=  \sum_{j,k=1}^{\infty}  \lim_{\alpha \rightarrow 0} a^2_{jk}(\alpha) = 0. \end{align}
The exchange in the order of limits and the sum in (\ref{Limit of HSnorm of A}) is permissable by the Lebesgue dominated convergence theorem since the summands satisfy $a^2_{jk}(\alpha) \leq 2 (c_{jk}^2(\alpha) + c^2_{jk})$ and
$\sum_{j,k=1}^{\infty}  a^2_{jk}(\alpha) < \infty$.
Now, since $\|B(\alpha) - B\|^2 \leq \|B(\alpha) - B\|^2_{HS} = \|A(\alpha)\|^2_{HS} \rightarrow 0$
as $\alpha \downarrow  0$, it follows that $B(\alpha)$ converges to $B$ in operator norm. \qed

\section{Tikhinov Regularized Canonical Correlation}

In the Tikhinov regularized approach to canonical correlation we replace the operators $\{ S_i^{1/2 \dag} \}_{i=1}^{2}$ with $\{ (S_i + \alpha I)^{-1/2} \}_{i=1}^{2}$ and then let
\begin{equation}\label{definition of Ralpha} R(\alpha) \equiv (S_1 + \alpha I)^{-1/2} S_{12} (S_2 + \alpha I)^{-1/2} \end{equation}
approximate the cross-correlation operator $R$ for $\alpha \in (0,a)$. Since the operators $\{ (S_1 + \alpha I)^{-1/2}, (S_2 + \alpha I)^{-1/2} \}$ are bounded and $S_{12}$ is Hilbert-Schmidt, it follows that $R(\alpha)$ is Hilbert-Schmidt. Now if we define $T(\alpha): \CH(K_2) \mapsto \CH(K_1)$ by
\begin{equation}\label{definition of Talpha} T(\alpha) \equiv \Gamma_1 R(\alpha) \Gamma^{-1}_2 = \Gamma_1 (S_1 + \alpha I)^{-1/2} S_{12} (S_2 + \alpha I)^{-1/2} \Gamma^{-1}_2 \end{equation} then as the regularization parameter $\alpha \downarrow 0$,
\begin{eqnarray}\label{Q alpha} R(\alpha) &=& (S_{1}+\alpha I_1)^{-1/2} S_{12} (S_{2}+\alpha I_2)^{-1/2}
\nonumber\\ &=& \sum_{j=1}^{\rank(S_1)} \sum_{k=1}^{\rank(S_2)} \frac{\gamma_{jk}}{\sqrt{(\lambda_{1j}+ \alpha)(\lambda_{2k}+ \alpha)}} \phi_{2k} \otimes_{L^2(E_2)} \phi_{1j}
\nonumber\\ &\rightarrow& \sum_{j=1}^{\rank(S_1)} \sum_{k=1}^{\rank(S_2)}  \frac{\gamma_{jk}}{\sqrt{\lambda_{1j} \lambda_{2k} }} \phi_{2k} \otimes_{L^2(E_2)} \phi_{1j}
\nonumber\\ &=& \Gamma^{-1}_1 T \Gamma_2. \end{eqnarray}
Therefore, by Therem \ref{HS Convergence theorem}, $R(\alpha)$ converges in operator norm
to the operator $\Gamma^{-1}_1 T \Gamma_1$ as $\alpha \downarrow 0$ and the continuity
of $\Gamma_1$ and $\Gamma_2$ ensures that
\beqs\label{T alpha} T(\alpha) = \Gamma_1 R(\alpha) \Gamma^{-1}_2 \rightarrow T \eeqs
as $\alpha \downarrow 0$, with convergence in terms of operator norm.

We will now show that the regularized canonical correlations along with the
regularized canonical variables converge to the canonical correlation and
variables defined from the Eubank and Hsing \citep{Eubank06}
methodology as $\alpha \downarrow 0$. In this regard, suppose that
$\{\rho_k(\alpha), f_{k}(\alpha), g_{k}(\alpha)\}_{k=1}^{\infty}$ is the singular system for
$R(\alpha)$ such that
\begin{equation*}\label{Q alpha2} R(\alpha) = \sum_{k=1}^{\infty} \rho_k(\alpha)
\left[ g_k(\alpha) \otimes_{\CH_2} f_k(\alpha) \right]. \end{equation*}
Then,
\begin{eqnarray*}\label{T alpha2} T(\alpha) = \Gamma_1 R(\alpha) \Gamma^{-1}_2  &=& \sum_{k=1}^{\infty} \rho_k(\alpha) \left[ (g_k(\alpha) \Gamma^{*}_2 ) \otimes_{\CH(K_2)} (\Gamma_1 f_{k}(\alpha)) \right]
\nonumber\\ &=& \sum_{k=1}^{\infty} \rho_k(\alpha) \left[ \tilde{g}_k(\alpha)  \otimes_{\CH(K_2)} \tilde{f}_k(\alpha) \right] \end{eqnarray*}
where $\tilde{f}_k(\alpha) = \Gamma_1 f_k(\alpha) \in \CH(K_1)$ and $\tilde{g}_k(\alpha) = \Gamma_2 g_k(\alpha) \in \CH(K_2)$.
Now by (\ref{RKHS Representation}) the canonical weight functions may be written as
\begin{equation*}\label{canonical weight formula} \tilde{f}_k(\alpha) = \sum_{j=1}^{\infty} \lambda_{1j} f_{kj}(\alpha) \phi_{1j}
~\text{ and }~ \tilde{g}_k(\alpha) = \sum_{j=1}^{\infty} \lambda_{2j} g_{kj}(\alpha) \phi_{2j} \end{equation*}
with
\begin{equation*}\label{canonical fourier coefficients} f_{kj}(\alpha) = \la \tilde{f}_k(\alpha),\phi_{1j} \ra_{\CH_1}
~\text{ and }~  g_{kj}(\alpha) = \la \tilde{g}_k(\alpha),\phi_{2j} \ra_{\CH_2}. \end{equation*}
Utilizing Theorem (\ref{Congruence Map Psi to L2X}) the corresponding regularized canonical variables in $L^2_{X_1}$ and $L^2_{X_2}$ are then
\begin{eqnarray*}\label{canonical variables formula} U_k(\alpha) &=& \Psi_1( \tilde{f}_k(\alpha) ) = \sum_{j=1}^{\rank(S_1)} f_{kj}(\alpha) \la X_1,\phi_{1j} \ra_{\CH_1} \text{ and } \nonumber\\ V_k(\alpha) &=& \Psi_2 ( \tilde{g}_k(\alpha) ) = \sum_{j=1}^{\rank(S_2)} g_{kj}(\alpha) \la X_2,\phi_{2j} \ra_{\CH_2}. \end{eqnarray*}

The continuity of the congruence mappings $\Psi_1$ and $\Psi_2$ ensures the convergence of the
regularized canonical variables $U_k(\alpha) = \Psi_1 ( \tilde{f}_k(\alpha) )$ and $V_k(\alpha) = \Psi_2 ( \tilde{g}_k(\alpha))$ to the true canonical variables, provided that the regularized canonical weight functions $\tilde{f}_k(\alpha) \in \CH(K_1)$ and
$\tilde{g}_k(\alpha) \in \CH(K_2)$ converge to the true canonical weight functions $\tilde{f}_k$ and $\tilde{g}_k$ as the regularization parameter tends to zero.
Thus, our tasks are to establish convergence of $\rho_k^2(\alpha)$ to $\rho^2_k$ and of the
the regularized RKHS functions $\{ \tilde{f}_k(\alpha),\tilde{g}_k(\alpha)\}$ to $\{\tilde{f}_k,\tilde{g}_k\}$
for all $k \geq 1$.  Concerning the convergence of the eigenvalues $\rho_k^2(\alpha)$ and the corresponding eigenprojection
operators we have the following result.

\begin{theorem}\label{CEGR THM1}
Let $\{\rho^2_k(\alpha),P_k(\alpha)\}$ and $\{\rho^2_k,P_k\}$ denote the eigenvalues and
corresponding eigenprojection operators for $T(\alpha)T^{*}(\alpha)$ and $TT^{*}.$ The following
limits hold as $\alpha \downarrow 0$
\begin{equation}\label{limit1}  0 \leq \rho_k^2(\alpha)  \uparrow \rho^2_k  \leq 1 \text{ as } \alpha \downarrow 0 ~ \text{ for all } k \geq 1 \text{ and } \end{equation}
\begin{equation}\label{limit2} \| P_k(\alpha) - P_k \| \leq  \| P_k(\alpha) - P_k \|_{HS} \rightarrow 0. \end{equation}
\end{theorem}

\noindent Proof: First note that $\rho^2_k(\alpha) < \rho^2_k < 1$ since $\|T(\alpha)\|^2 < \|T\|^2 \leq 1$. (see Proposition A.3 of
Eubank and Hsing \citep{Eubank06}). Now to see that (\ref{limit1}) holds, fix $k \geq 1$. Since $\rho_k^2(\alpha)$ and $\rho^2_k$ are the
$k^{th}$ eigenvalues for $T(\alpha)T^{*}(\alpha)$ and $TT^{*}$
\begin{equation*}\label{proof of eigenvalue convergence} \left( \rho^2_k - \rho_k^2(\alpha) \right) = | \rho^2_k - \rho_k^2(\alpha) |
\leq \| TT^{*} - T(\alpha)T^{*}(\alpha) \| \downarrow 0 \end{equation*}
as $\alpha \downarrow 0$.  In order to show that $P_k(\alpha) \rightarrow P_k$ in operator norm, let $\Gamma_{r,k}$ be a circle centered at $\rho^2_k$
with radius $r$ chosen so that $\Gamma_{r,k}$ encloses $\rho^2_k$ and no other eigenvalues of $TT^{*}$.
Suppose that $\{R(\alpha,z), R(z)\}$ are the resolvents of $\{T(\alpha)T^{*}(\alpha),TT^{*}\}$, respectively.
Since $\|T(\alpha)T^{*}(\alpha) - TT^{*}\| \rightarrow 0$ as $\alpha \downarrow 0$, it follows from
Theorem \ref{continuity of specturm} in the appendix
that there exists $\alpha_{0} >0$ such that whenever $0 < \alpha < \alpha_{0}$, $\Gamma_{r,k}$ encloses
$\rho^2_k(\alpha)$ and no other eigenvalues of $T(\alpha)T^{*}(\alpha)$.
Furthermore, for any  $\epsilon > 0$ we may take $\alpha_{0}$ to be sufficiently small to ensure that
$\|T(\alpha)T^{*}(\alpha) - TT^{*}\| < \epsilon$. Relation (\ref{Proj2}) from the appendix then has the consequence that
\begin{equation*}\label{proof that projection operators converge} \| P_k(\alpha) - P_k \|_{HS} \leq r \sup_{z \in \Gamma_{r,k}}
\left\{ \frac{\|T(\alpha)T^{*}(\alpha) - TT^{*} \|_{HS}  \|R(z)\|^2_{HS} }
{ 1- \|T(\alpha)T^{*}(\alpha) - TT^{*} \|_{HS}  \|R(z)\|_{HS}} \right\}. \end{equation*}
Thus, if $M(r,k) \equiv \sup_{z \in \Gamma_{r,k}} \|R(z)\|_{HS}$ and $\epsilon >0$ are chosen
so that $\epsilon < \frac{1}{2 M(r,k)}$, $\| P_k(\alpha) - P_k \|_{HS} \leq 2 r M^2(r,k) \epsilon$
and hence $\| P_k(\alpha) - P_k \| < \| P_k(\alpha) - P_k \|_{HS} \rightarrow 0$ as $\alpha \downarrow 0$. \qed

It remains to show that for $k \geq 1$,
$\tilde{f}_k(\alpha) = \Gamma_1(f_k(\alpha))$ and $\tilde{g}_k(\alpha) = \Gamma_2(g(\alpha))$
approach $\tilde{f}_k \in \CH(K_1)$ and $\tilde{g}_k \in \CH(K_2)$
from the singular system $\{\rho_k,\tilde{f}_k,\tilde{g}_k \}$ of $T$. We note, however,
that eigenvectors associated with any operator are not defined uniquely.
For example, if $\theta$ is an eigenvector for an arbitrary self-adjoint operator $A$,
then $-\theta$ is also an eigenvector.  In order to properly establish what we mean by convergence assume, WLOG, that for all $\alpha >0$,
$\tilde{f}_k(\alpha)$ be chosen so that $ \la \tilde{f}_k(\alpha),\tilde{f}_k \ra_{\CH(K_1)} \geq 0$, with a similar convention applied to $\tilde{g}_k(\alpha)$.
The theorem below concerns convergence in the case that the eigenspaces associated with
$\{\tilde{f}_{k},\tilde{g}_{k}\}$ are 1--dimensional.  Subsequently, we will discuss the higher dimensional case.

\begin{theorem}\label{CEGR THM2}
Assume that the eigenspaces associated with the eigenvectors $\tilde{f}_{k}$ and $\tilde{g}_{k}$
are one dimensional with $k \in \BBZ$. Then, as $\alpha \downarrow 0$
\begin{equation*}\label{rspcc weight functions converge} \|\tilde{f}_k(\alpha) - \tilde{f}_k \|_{\CH(K_1)} \rightarrow 0
~\text{ and }~ \| \tilde{g}_k(\alpha) - \tilde{g}_k \|_{\CH(K_2)} \rightarrow 0.\end{equation*}
\end{theorem}

\noindent Proof: For fixed $k \in \BBZ$, it suffices to show that
$\| \tilde{f}_k(\alpha) - \tilde{f}_k\|_{\CH(K_1)} \rightarrow 0$.
Since the eigenspaces are one-dimensional it follows that $P_k(\alpha) =
\left[ \tilde{f}_k(\alpha) \otimes_{\CH(K_1)} \tilde{f}_k(\alpha) \right]$ and
$P_k = \left[ \tilde{f}_k \otimes_{\CH(K_1)} \tilde{f}_k \right]$. Now, notice that
\begin{align}\label{relationship between projections and eigenvectors}
\| P_k(\alpha) - P_k \|_{HS} &= \la P_k(\alpha) - P_k, P_k(\alpha) - P_k \ra_{HS}
\notag\\ &= 2 - 2 \la P_k(\alpha),P_k \ra_{HS}
\notag\\ &= 2 - 2 \la \tilde{f}_k(\alpha) \otimes_{\CH(K_1)} \tilde{f}_k(\alpha),\tilde{f}_k \otimes_{\CH(K_1)} \tilde{f}_k \ra_{HS}
\notag\\ &= 2 - 2 \la \tilde{f}_k(\alpha),\tilde{f}_k \ra^2_{\CH(K_1)}
\notag\\ &= 2(1- \la \tilde{f}_k(\alpha),\tilde{f}_k \ra_{\CH(K_1)})(1+ \la \tilde{f}_k(\alpha),\tilde{f}_k \ra_{\CH(K_1)})
\notag\\ &= \|\tilde{f}_k(\alpha) - \tilde{f}_k\|^2_{\CH(K_1)}(1+ \la \tilde{f}_k(\alpha),\tilde{f}_k \ra_{\CH(K_1)}).
\end{align}
Furthermore, as $\la \tilde{f}_k(\alpha),\tilde{f}_k \ra_{\CH(K_1)} \geq 0 \longrightarrow (1+ \la \tilde{f}_k(\alpha),\tilde{f}_k \ra_{\CH(K_1)}) \geq 1$, hence
\begin{equation*}\label{norm of eigenvectors less that projection ops}
\|\tilde{f}_k(\alpha) - \tilde{f}_k\|^2_{\CH(K_1)} \leq \| P_k(\alpha) - P_k \|_{HS}. \end{equation*}
Since $\| P_k(\alpha) - P_k \|_{HS} \rightarrow 0$ as $\alpha \downarrow 0$, it follows
that $\|\tilde{f}_k(\alpha) - \tilde{f}_k\|^2_{\CH(K_1)} \rightarrow 0$. \qed

It should be noted that if $\la \tilde{f}_k(\alpha),\tilde{f}_k \ra_{\CH(K_1)} \leq 0$ instead,
then $(1 - \la \tilde{f}_k(\alpha),\tilde{f}_k \ra_{\CH(K_1)}) \geq 1$ and from (\ref{relationship between projections and eigenvectors})
we would have
\begin{eqnarray*}\label{alternative if convention is not accepted}
\| P_k(\alpha) - P_k \|_{HS} &=& 2(1 + \la \tilde{f}_k(\alpha),-\tilde{f}_k \ra_{\CH(K_1)})(1 - \la \tilde{f}_k(\alpha),-\tilde{f}_k \ra_{\CH(K_1)})
\nonumber\\ &=& (1 - \la \tilde{f}_k(\alpha),\tilde{f}_k \ra_{\CH(K_1)}) \|\tilde{f}_k(\alpha) - (- \tilde{f}_k) \|^2_{\CH(K_1)}.
\end{eqnarray*}
Hence,
\begin{equation*}\label{norm of eigenvectors less projection ops 2 case}
\|\tilde{f}_k(\alpha) - (-\tilde{f}_k) \|^2_{\CH(K_1)} \leq \| P_k(\alpha) - P_k \|_{HS}
\rightarrow 0 ~\text{ as }~ \alpha \downarrow 0 \end{equation*}
and $\tilde{f}_k(\alpha)$ would converge to $(-\tilde{f}_k)$ instead.

When the eigenspaces have dimension larger than 1, it is possible to find infinitely many eigenspace invariant rotations $\Theta \in \CB(\CH(K_1))$ so that $\tilde{f}'_{k} = \Theta \tilde{f}_{k}$ is still an eigenvector of $TT^{*}$ with eigenvalue $\rho^2_k$, yet $\|\tilde{f}_k(\alpha) - \tilde{f}'_k \| \nrightarrow 0$ as
$\alpha \downarrow 0$ (see Kato \cite{Kato80} p. 98-100).

Theorems \ref{CEGR THM1} and \ref{CEGR THM2} ensure that when $T$ is simple, the singular system $\{\rho_k(\alpha),\tilde{f}_k(\alpha),\tilde{g}_k(\alpha)\}$
of $T(\alpha)$ converges to the singular system $\{\rho_k,\tilde{f}_k,\tilde{g}_k\}$ of $T$ as the regularization parameter $\alpha \downarrow 0$. This is a positive development provided the singular value decomposition of $T(\alpha)$ can be estimated.  However, the singular value decomposition of $T(\alpha)$ entails the eigenvalue-eigenvector decomposition of e.g., the operator
\begin{eqnarray}\label{exact operator for TTstar} \CT_1(\alpha) &\equiv& T(\alpha)T^{*}(\alpha) = \Gamma_1 R(\alpha) R^{*}(\alpha) \Gamma^{-1}_1
\notag\\ &=& \Gamma_1 (S_{1}+\alpha I)^{-1/2} S_{12} (S_{2}+\alpha I)^{-1} S_{21} (S_{1}+\alpha I)^{-1/2} \Gamma^{-1}_1.\end{eqnarray}

Since, $\Gamma_1$ is unknown in (\ref{exact operator for TTstar}) we might estimate it using
\begin{equation*}\label{possible estimator for gamma} \hat{\Gamma}_{1n}(m) = \sum_{i=1}^{m} \sqrt{\hat{\lambda}_{1in}} \hat{P}_{1in} \end{equation*}
with $\{\hat{\lambda}_{1in},\hat{P}_{1in}\},$ the estimated eigenvalues and corresponding eigenprojection operators for $\hat{S}_{1n}$ and $m$ some integer. This raises the question of how to select $m$ and, for large $m$, $ \hat{\Gamma}_{1n}(m)$ is approximately $\hat{S}_{1n}^{1/2}$ whose compact nature is what prompted us to regularize from the beginning.

Since we are already utilizing Tikhinov regularization, a possible remedy for our problem is
to replace $\Gamma_1$ with $(S_1 + \alpha I)^{1/2}$. This produces the operator
\begin{equation}\label{Definition for CS1} \CS_1(\alpha) \equiv  S_{12} (S_{2}+\alpha I)^{-1} S_{21} (S_{1}+\alpha I)^{-1}\end{equation}
whose domain is $\ker(S_1)^{\perp}$ rather than $\CH(K_1)$. One advantage of $\CS_1(\alpha)$ is that
$\Ima(\CS_1(\alpha)) \subseteq \Ima(S_1) \subseteq F_1$ and hence the eigenfunctions of $\CS_1(\alpha)$ satisfy the Picard criteria. To see this, note by the infinite dimensional extension of the result from Khatri \citep{Khatri76} we have that $S_1 S_1^{\dag} S_{12} = S_{12}$ and hence  $\CS_1(\alpha) =S_1 S_1^{\dag} \CS_1(\alpha)$ (see King \citep{King10}).  Note that the operator $\CS_1(\alpha)$ is self-adjoint since
\begin{eqnarray*}\label{Proof that CS1 is HS} \CS_1(\alpha) &=& \sum_{j=1}^{\infty} \sum_{k=1}^{\infty} \frac{\gamma^2_{jk}}{(\lambda_{1j}+\alpha)(\lambda_{2k}+\alpha)} \left[ \phi_{1j} \otimes_{\CH_1} \phi_{1k} \right]
\nonumber\\ &=& (S_1+\alpha I)^{-1} S_{12} (S_2+\alpha I)^{-1} S_{21} = \CS^{*}_1(\alpha).  \end{eqnarray*}
Furthermore, as $S_{12}$ is a factor in $\CS_1(\alpha)$, the operator is Hilbert-Schmidt and
hence admits an eigenvalue-eigenvector decomposition
\begin{equation*}\label{eigenfunction expansion for CS1} \CS_1(\alpha) = \sum_{j=1}^{\infty}
{\rho}^{2}_j(\alpha) \left[ f_j(\alpha) \otimes_{\CH_1} f_j(\alpha) \right]. \end{equation*}
Now the question becomes how can an operator whose domain and range are subsets of $\CH_1$,
approximate an operator whose domain and range are subsets of $\CH(K_1)$.  The answer to this question was fundamentally answered by Nasheed and Wahba \citep{Wahba74} when it was proved that the collection of functions in $\CH(K_1)$ is the same as $\overline{\Ima(S_1)^{1/2}}$, except with alternate norm and inner product.  As the collection of eigenfunctions $\{ f_j(\alpha) \}_{j=1}^{\infty}$ reside in $\Ima(S_1)^{1/2}$, they also have ``dual citizenship'' in $\CH(K_1)$. We may therefore regard the eigenfunction sequence $\{ f_j(\alpha) \}$ as residing in $\CH(K_1)$, provided that we norm the eigenfunctions correctly.  If we treat the eigenfuctions $f_j(\alpha)$ as citizens of $\CH(K_1)$, for notational consistency we will denote them by $\tilde{f}_j(\alpha)$ with $\{f_j(\alpha) = \tilde{f}_j(\alpha) \}_{j=1}^{\rank(\CS_1(\alpha))}$.  There are therefore two possible views one
may adopt concerning the operator $\CS_1(\alpha)$:
\begin{enumerate}
\item In the first view of $\CS_1(\alpha)$, we treat the operator as a self-adjoint mapping
in $\ker(S_1)^{\perp} \subseteq L^2(E_1)$.
\item In the second viewpoint, the operator is treated as a self-adjoint mapping on $\CH(K_1)$ with
$\CS_1(\alpha)$ regarded as ``two perturbations'' distant from the the operator $TT^{*}$, which
is its ultimate intended target of approximation.
\end{enumerate}

\noindent When the second viewpoint for $\CS_1(\alpha)$ is adopted, the operator $\CS_1(\alpha)$ is representable by the $\CH(K_1)$ based operator
\begin{equation*}\label{slick form of the RKHS representation of CS1} \CS_1(\alpha) = \Gamma_1 S_1^{1/2 \dag} S_{12} (S_{2}+\alpha I)^{-1} S_{21} (S_{1}+\alpha I)^{-1}S_1^{1/2} \Gamma_1^{-1}.\end{equation*}
We then see that as $\alpha \downarrow 0$,
\begin{eqnarray*}\label{proof after assignment} \Gamma_1^{-1} \CS_1(\alpha) \Gamma_1 &=& S_1^{1/2 \dag} S_{12} (S_{2}+\alpha I)^{-1} S_{21} (S_{1}+\alpha I)^{-1} S_1^{1/2}
\nonumber\\ &=&  \sum_{j=1}^{\infty} \sum_{k=1}^{\infty} \frac{\gamma^2_{jk}}{(\lambda_{1j}+ \alpha)(\lambda_{2k}+ \alpha)}
\left[ \phi_{1j} \otimes_{\CH_1} \phi_{1k} \right]
\nonumber\\ &\rightarrow&  \sum_{j=1}^{\infty} \sum_{k=1}^{\infty} \frac{\gamma^2_{jk}}{(\lambda_{1j})(\lambda_{2k})}
\left[ \phi_{1j} \otimes_{\CH_1} \phi_{1k} \right]
\nonumber\\ &=&  \Gamma_1^{-1} T T^{*} \Gamma_1. \end{eqnarray*}
Therefore, since $\CS_1(\alpha)$ is Hilbert-Schmidt, Theorem \ref{HS Convergence theorem} ensures that as the regularization
parameter $\alpha \downarrow 0$,
\begin{equation*}\label{cs convergence in norm} \| \Gamma_1^{-1} \CS_1(\alpha) \Gamma_1 - \Gamma_1^{-1} TT^{*} \Gamma_1\|_{\CH_1} = \| \CS_1(\alpha)- TT^{*}\|_{\CH(K_1)} \rightarrow 0. \end{equation*}

\section{Asymptotic Properties for Tikhinov Regularization}

In this section we will consider the asymptotics associated with the sample
estimators of the operators $\CS_1(\alpha)$. The asymptotics associated with the
operator $\CS_1(\alpha)$ rely heavily on perturbation theory concepts discussed in Dauxois et al. \citep{DNR82} as well as delta method theory for random operators discussed in Cupidon et al. \citep{Cupidon07}.

To begin, we suppose that a random sample $X_1, X_2, \ldots X_n$ of independent, identically distributed copies of $X \in L^2(E)$ are observed. The sample estimator associated with the covariance operator of $X$ is given by
\begin{equation}\label{S hat Defn tikhinov section} \hat{S}_{n} = \frac{1}{n} \sum_{i=1}^{n} (X_i - \bar{X}) \otimes_{\CH} (X_i - \bar{X})\end{equation}
and the continuous mapping theorem along with the law of large numbers ensures that
$\hat{S}_{in}  = \Upsilon_i \hat{S}_{n} \Upsilon_i \cip S_{i}$
for $i=1,2$ and $\hat{S}_{12n}  = \Upsilon_1 \hat{S}_{n} \Upsilon_2 = \hat{S}^{*}_{21n} \cip S_{12}$ as $n \rightarrow \infty$. For Tikhinov
regularization we will have need of the function $\varphi_{\alpha}(z) \equiv (z+\alpha)^{-1},$ which is analytic for all points in the complex plane,
except for a pole at $z = - \alpha.$
Consequently, the disk $D = \{ z \in \BBC | \min_{0 \leq x \leq \|S\|} |z - x| < \frac{\alpha}{2} \}$
contains the spectra of $S$ and the function $\varphi_{\alpha}$ is analytic on $D$.
It follows by the continuous mapping theorem that $\varphi_{\alpha}(\hat{S}_{in}) \cip \varphi_{\alpha}(S_i)$ as $n \rightarrow \infty$ for $i = 1,2$.

As a consequence of the continuous mapping theorem and the central limit theorem for Hilbert space operators (see Dauxois et al. \citep{DNR82})  we have that
\beq\label{conv in dist tikhinov section} \sqrt{n}( \Upsilon_i \hat{S}_{n} \Upsilon_j - \Upsilon_i S \Upsilon_j) = \sqrt{n}( \hat{S}_{ijn} - S_{ij} ) \cid  \Upsilon_i \CN \Upsilon_j = \CN_{ij} \eeq
where, for $i,j = 1,2$, $\CN_{ij} \in \CK_{HS}(\CH_i,\CH_j)$ is a Gaussian random operator that has mean zero and variance
\begin{equation}\label{variance tikhinov section} \Sigma_{ij} \equiv {\rm E}\{ (X_i \otimes_{\CH_i} X_i - S_i) \otimes_{HS} (X_j \otimes_{\CH_j} X_j - S_j) \} \end{equation}
and $\CN_{ii} \equiv \CN_{i}.$
Furthermore, by the delta method result from Cupidon et al. \citep{Cupidon07} it follows that for $i = 1,2$
\beq\label{tikhinov function in tikhinov section} \sqrt{n} \left\{ \varphi_{\alpha}( \hat{S}_{in} )  - \varphi_{\alpha}(S_i) \right\} \cid \varphi'_{\alpha}( \CN_i) \eeq
where the limit in $\CK_{HS}(\CH_i)$ has zero mean and is distributed as
\begin{eqnarray}\label{delta method tikhinov section} \varphi'_{\alpha}(\CN_i) = &-& \sum_{k =1}^{\infty} (\lambda_{ik} + \alpha)^{-2} P_{ik} \CN_i P_{ik} \nonumber\\
&+&  \sum_{j \neq k} \frac{1}{(\lambda_{ik} + \alpha)(\lambda_{ij} + \alpha)} P_{ij} \CN_i P_{ik} \end{eqnarray}
with $\{\lambda_{ik},P_{ik}\}_{k=1}^{\infty}$, the eigenvalues and eigenprojection operators corresponding to $S_i$, $i = 1,2$ (see Appendix).

The sample version of the operator $\CS_1(\alpha)$ is then defined by
\begin{equation}\label{CEGR sample estimator for CS} \hat{\CS}_{1n}(\alpha) \equiv \hat{S}_{12n} (\hat{S}_{2n} + \alpha I)^{-1} \hat{S}_{21n} (\hat{S}_{1n} + \alpha I)^{-1}.\end{equation}
The asymptotic analysis of $\sqrt{n}(\hat{\CS}_{1n}(\alpha) - \CS_{1}(\alpha))$ follows from a product rule application of the delta method similar to that in Cupidon, et al. \citep{Cupidon07}. In this regard we introduce the following Gaussian elements in the set $\CK_{HS}(\CH(K_1))$ of Hilbert-Schmidt operators on $\CH(K_1)$
\begin{eqnarray} \CG_{11}(\alpha) &\equiv& \CN_{12} \varphi_{\alpha}(S_2) S_{21} \varphi_{\alpha}(S_1),\nonumber\\
\CG_{12}(\alpha) &\equiv& S_{12} \varphi'_{\alpha}(\CN_2) S_{21} \varphi_{\alpha}(S_1),\nonumber\\
\CG_{13}(\alpha) &\equiv& S_{12} \varphi_{\alpha}(S_2) \CN_{21} \varphi_{\alpha}(S_1),\nonumber\\
\CG_{14}(\alpha) &\equiv& S_{12} \varphi_{\alpha}(S_2) S_{21} \varphi'_{\alpha}(\CN_1),\nonumber\\
\CG_1(\alpha) &\equiv& \sum_{k=1}^{4} \CG_{1k}(\alpha).\label{gaussian elemnts for tikhinov asymptotics}\end{eqnarray}

\begin{corollary}\label{complex theorem in tikhinov section}
If ${\rm E} \| X \|_{L^2(E)}^4 < \infty$, then as $n \rightarrow \infty$
\beq\label{complex distribution tikhinov section} \sqrt{n}(\hat{\CS}_{1n}(\alpha) - \CS_{1}(\alpha)) \cid \CG_1(\alpha).\eeq
\end{corollary}

\noindent Proof: Define
\begin{eqnarray*} \hat{\CA}_{11}(\alpha) &\equiv& \left[ \hat{S}_{12n} - S_{12} \right] \varphi_{\alpha}(\hat{S}_{2n}) \hat{S}_{21n} \varphi_{\alpha}(\hat{S}_{1n}),\nonumber\\
 \hat{\CA}_{12}(\alpha) &\equiv& S_{12} \left[ \varphi_{\alpha}(\hat{S}_{2n}) - \varphi_{\alpha}(S_{2}) \right] \hat{S}_{21n} \varphi_{\alpha}(\hat{S}_{1n}),\nonumber\\
 \hat{\CA}_{13}(\alpha) &\equiv& S_{12} \varphi_{\alpha}(S_{2}) \left[ \hat{S}_{21n} - S_{21} \right]  \varphi_{\alpha}(\hat{S}_{1n}),\nonumber\\
 \hat{\CA}_{14}(\alpha) &\equiv& S_{12} \varphi_{\alpha}(S_{2}) S_{21} \left[ \varphi_{\alpha}(\hat{S}_{1n}) - \varphi_{\alpha}(S_{1}) \right].
\label{1st differences in proof for S Tikhinov asymptotics}\end{eqnarray*}
Notice that the difference $\sqrt{n}(\hat{\CS}_{1n}(\alpha) - \CS_{1}(\alpha))$ can be expanded so that
\begin{equation*}\label{expansion of difference of cs in Tikhinov proof} \sqrt{n}(\hat{\CS}_{1n}(\alpha) - \CS_{1}(\alpha)) =
\sqrt{n} \left[ \sum_{j=1}^{4} \hat{\CA}_{1j}(\alpha) \right]. \end{equation*}
The application of (\ref{conv in dist tikhinov section}), (\ref{tikhinov function in tikhinov section}) and Slutsky's Theorem then ensure that
\begin{equation*} \sqrt{n} \left[ \sum_{j=1}^{4} \hat{\CA}_{1j}(\alpha) \right] \cid \CG_1(\alpha)\end{equation*}
since, for example, the term $\hat{\CA}_{11}(\alpha)$ consists of the factor $\sqrt{n} \left[ \hat{S}_{12n} - S_{12} \right] \cid \CN_{12}$
right-multiplied by the factor
\begin{equation*} \varphi_{\alpha}(\hat{S}_{2n}) \hat{S}_{21n} \varphi_{\alpha}(\hat{S}_{1n}) \cip \varphi_{\alpha}(S_{2}) S_{21} \varphi_{\alpha}(S_{1}).\end{equation*} \qed

As a result of Corollary \ref{complex theorem in tikhinov section}, we see that $ \hat{\CS}_{1n}(\alpha)$  is a consistent estimator of $\CS_{1}(\alpha)$ as
\begin{equation}\label{consistency tikhinov section} \| \hat{\CS}_{1n}(\alpha) - \CS_{1}(\alpha) \| = \CO_{P}(n^{-1/2}) \cip 0.\end{equation}
However, note that as long as the regularization parameter $\alpha > 0$, $\| \hat{\CS}_{1n}(\alpha) - T T^{*} \| \ncip 0.$
In fact, by the triangle inequality we have
\begin{equation}\label{consistency triangle inequality tikhinov section}  \| \hat{\CS}_{1n}(\alpha) - T T^{*} \| \leq \| \hat{\CS}_{1n}(\alpha) - \CS_{1}(\alpha) \| + \| \CS_{1}(\alpha) - T T^{*} \|. \end{equation}
The first term on the right-hand side of (\ref{consistency triangle inequality tikhinov section}) can be viewed as a random error that originates
from using a sample estimator of $\CS_{1}(\alpha)$. This term tends to zero as $n \rightarrow \infty$
by (\ref{consistency tikhinov section}).  On the other hand, the second term on the right hand side of
(\ref{consistency triangle inequality tikhinov section}) is a deterministic error that arises from the
regularized approximation of $T T^{*}$. This latter term will only become negligible if $\alpha \downarrow 0$.

Since the limiting distribution for $\sqrt{n}(\hat{\CS}_{1n}(\alpha) - \CS_{1}(\alpha))$ has been established, we may establish the limiting distributions associated with sample estimators for the $k^{th}$ regularized canonical correlation and associated projection operator and weight functions.  The quantities
of interest are $\sqrt{n} \left\{ \hat{\rho}_{kn}(\alpha) - \rho_{k}(\alpha) \right\}$,
$\sqrt{n} \left\{ \hat{\tilde{P}}_{1kn}(\alpha) - \tilde{P}_{1k}(\alpha) \right\}$
and $\sqrt{n} \left\{ \hat{\tilde{f}}_{kn}(\alpha) - \tilde{f}_{k}(\alpha) \right\}$ where
$\{\rho_{k}(\alpha),\tilde{P}_{1k}(\alpha),\tilde{f}_{k}(\alpha)\}$ denote the eigenvalues, eigenprojections and eigenvectors for $\CS_{1}(\alpha)$
and $\{\hat{\rho}_{kn}(\alpha),\hat{\tilde{P}}_{1kn}(\alpha),\hat{\tilde{f}}_{kn}(\alpha)\}$ denote the same for $\hat{\CS}_{1n}(\alpha)$.

\begin{theorem}\label{theorem for projections and eigenvectors}
Suppose that ${\rm E} \| X \|_{L^2(E)}^4 < \infty$. Then, as $n \rightarrow \infty$
\beq\label{limiting distribution of P in Tikhinov section} \sqrt{n} \left\{ \hat{\tilde{P}}_{1kn}(\alpha) - \tilde{P}_{1k}(\alpha) \right\}
\cid  \tilde{P}_{1k}(\alpha)  \CG_1(\alpha) \tilde{Q}_{1k}(\alpha)  + \tilde{Q}_{1k}(\alpha)  \CG_1(\alpha) \tilde{P}_{1k}(\alpha) \eeq
where $\CG_1(\alpha)$ is as in (\ref{gaussian elemnts for tikhinov asymptotics}) and
\beqs\label{definition of Qk alpha}  \tilde{Q}_{1k}(\alpha) = \sum_{j \neq k} \frac{1}{\rho_{j}(\alpha) - \rho_{k}(\alpha)} \tilde{P}_{1j}(\alpha). \eeqs
In the case that $\rank(\tilde{P}_{1k}(\alpha)) = 1$,
\beq\label{limiting distribution of eigenfunctions} \sqrt{n} \left\{ \hat{\tilde{f}}_{kn}(\alpha) - \tilde{f}_{k}(\alpha) \right\}
\cid  \tilde{Q}_{1k}(\alpha)  \CG_1(\alpha) \tilde{f}_{k}(\alpha). \eeq
\end{theorem}

\noindent Proof: For each $k \in \BBZ$, let $\Gamma_k$ denote a circle that encloses the eigenvalue $\rho_{k}(\alpha)$ but no other eigenvalue eigenvalue of $\CS_{1}(\alpha)$.  It follows from developments in the appendix that
\begin{equation}\label{tikhinov asymptotics on projections ops in proof} \sqrt{n} \left\{ \hat{\tilde{P}}_{1kn}(\alpha) - \tilde{P}_{1k}(\alpha) \right\} =  \frac{\sqrt{n}}{2 \pi i} \oint_{\Gamma_k}  R(z)(\hat{\CS}_{1n}(\alpha) - \CS_{1}(\alpha)) R(z) H_n(z,\alpha) dz \end{equation}
with $H_n(z,\alpha) \equiv \sum_{j=0}^{\infty} \left\{ (\CS_{1}(\alpha) - \hat{\CS}_{1n}(\alpha) )R(z) \right\}^{j}$ and $R(z)$ the resolvent
of $\CS_{1}(\alpha)$.  Now since the integrand in (\ref{tikhinov asymptotics on projections ops in proof}) can be expanded into
$$ R(z)(\hat{\CS}_{1n}(\alpha) - \CS_{1}(\alpha)) R(z) H_n(z,\alpha) = R(z)(\hat{\CS}_{1n}(\alpha) - \CS_{1}(\alpha)) R(z) + M(z,\alpha)$$
with
\begin{align*}\label{Mn} M(z,\alpha) &\equiv R(z)(\hat{\CS}_{1n}(\alpha) - \CS_{1}(\alpha)) R(z) \sum_{j=1}^{\infty} \left\{(\hat{\CS}_{1n}(\alpha) - \CS_{1}(\alpha)) R(z)\right\}^{j}
\notag\\ &= R(z)(\hat{\CS}_{1n}(\alpha) - \CS_{1}(\alpha)) R(z)(\hat{\CS}_{1n}(\alpha) - \CS_{1}(\alpha)) R(z) \notag\\ &+R(z)(\hat{\CS}_{1n}(\alpha) - \CS_{1}(\alpha)) R(z)(\hat{\CS}_{1n}(\alpha) - \CS_{1}(\alpha)) R(z)(\hat{\CS}_{1n}(\alpha) - \CS_{1}(\alpha)) R(z)+ \cdots
\notag\\ &= \CO_P(n^{-1}). \end{align*}
It follows that
\begin{equation}\label{Asym 2} \sqrt{n} \left\{ \hat{\tilde{P}}_{1kn}(\alpha) - \tilde{P}_{1k}(\alpha) \right\} =  \frac{\sqrt{n}}{2 \pi i} \oint_{\Gamma_k}  R(z)(\hat{\CS}_{1n}(\alpha) - \CS_{1}(\alpha)) R(z) dz + \CO_P(n^{-1/2}).\end{equation}
We may now focus attention on the lead term in (\ref{Asym 2}).  From Corollary \ref{complex theorem in tikhinov section} and the continuous mapping theorem it follows that
\begin{equation}\label{Asym 3} \sqrt{n} \left\{ \hat{\tilde{P}}_{1kn}(\alpha) - \tilde{P}_{1k}(\alpha) \right\} \cid  \frac{1}{2 \pi i} \oint_{\Gamma_k}  R(z) \CG_1(\alpha) R(z) dz.\end{equation}
To simplify the last expression we write
\begin{equation*}\label{Resolvent expansion} R(z) = \sum_{k=1}^{\infty} \frac{1}{\rho_{k}(\alpha) - z} \tilde{P}_{1k}(\alpha) + \CO( (\rho_{k}(\alpha) -z)^{-2})\end{equation*}
and all but the lead term will vanish when the contour integral is taken due to (\ref{Proof showing nilpotent terms can be ignored}).
The integrand in (\ref{Asym 3}) can then be simplified as
\begin{equation}\label{Sandwich1} R(z) \CG_1(\alpha) R(z)  = \sum_{k=1}^{\infty} \sum_{j=1}^{\infty}  \frac{1}{(\rho_{k}(\alpha) - z)(\rho_{j}(\alpha) - z)}\tilde{P}_{1k}(\alpha) \CG_1(\alpha) \tilde{P}_{1j}(\alpha). \end{equation}
Applying the Cauchy integral formula to (\ref{Sandwich1}) ensures that
\begin{equation*}\label{integral1} \frac{1}{2 \pi i} \oint_{\Gamma_k}  \frac{dz}{(\rho_{k}(\alpha) - z)(\rho_{j}(\alpha) - z)} = \frac{1}{2 \pi i} \oint_{\Gamma_k}  \frac{(\rho_{i}(\alpha) - z)dz}{(\rho_{k}(\alpha) - z)(\rho_{j}(\alpha) - z)(\rho_{i}(\alpha) - z)} \end{equation*}
and the only case where the integral is non-zero is when exactly one of $\rho_{k}(\alpha)$ or $\rho_{j}(\alpha)$ is not equal to $\rho_{i}(\alpha)$. When, for example, $\rho_{k}(\alpha) = \rho_{i}(\alpha)$ and $\rho_{i}(\alpha) \neq \rho_{j}(\alpha)$ we have
\begin{equation*}\label{integral 2} \frac{1}{2 \pi i} \oint_{\Gamma_k}  \frac{(\rho_{i}(\alpha) - z)}{(\rho_{k}(\alpha) - z)(\rho_{j}(\alpha) - z)(\rho_{i}(\alpha) - z)} dz = \frac{1}{(\rho_{j}(\alpha) - \rho_{k}(\alpha))} \end{equation*}
and hence
\begin{eqnarray*}\label{Asym 4}  \sqrt{n} \left\{ \hat{\tilde{P}}_{1kn}(\alpha) - \tilde{P}_{1k}(\alpha) \right\} &\cid&
\sum_{i=1}^{\infty} \sum_{j \neq k} \frac{\delta_{ik}}{(\rho_{j}(\alpha) - \rho_{k}(\alpha))} \tilde{P}_{1i}(\alpha) \CG_1(\alpha) \tilde{P}_{1j}(\alpha)  \nonumber\\ &+& \sum_{j=1}^{\infty} \sum_{i \neq k} \frac{\delta_{jk}}{(\rho_{i}(\alpha) - \rho_{k}(\alpha))} \tilde{P}_{1i}(\alpha) \CG_1(\alpha) \tilde{P}_{1j}(\alpha)  \nonumber\\ &=& \tilde{P}_{1k}(\alpha)  \CG_1(\alpha) \tilde{Q}_{1k}(\alpha)  + \tilde{Q}_{1k}(\alpha)  \CG_1(\alpha) \tilde{P}_{1k}(\alpha)\end{eqnarray*}
which establishes (\ref{limiting distribution of P in Tikhinov section}).

To obtain the limiting distribution of $\sqrt{n} \left\{ \hat{\tilde{f}}_{kn}(\alpha) - \tilde{f}_{k}(\alpha) \right\},$
first observe that an application of Theorem \ref{continuity of specturm}
ensures that for large $n$ and probability tending to $1$,  $\rank(\hat{\tilde{P}}_{1kn}(\alpha)) = 1$.
Thus, we may write $\hat{\tilde{P}}_{1kn}(\alpha) = \left[ \hat{\tilde{f}}_{kn}(\alpha) \otimes_{\CH(K_1)} \hat{\tilde{f}}_{kn}(\alpha) \right]$ and hence
\begin{eqnarray*}\label{Proj op and eigenvector relation} && \la \hat{\tilde{P}}_{1kn}(\alpha) - \tilde{P}_{1k}(\alpha), \tilde{P}_{1k}(\alpha) \ra_{HS} =  \la \hat{\tilde{P}}_{1kn}(\alpha), \tilde{P}_{1k}(\alpha) \ra_{HS} - 1 \nonumber\\
&=&  \la ( \hat{\tilde{f}}_{kn}(\alpha) \otimes_{\CH(K_1)} \hat{\tilde{f}}_{kn}(\alpha) ), (\tilde{f}_{k}(\alpha) \otimes_{\CH(K_1)} \tilde{f}_{k}(\alpha) ) \ra_{HS} - 1 \nonumber\\
&=&  \la \hat{\tilde{f}}_{kn}(\alpha), \tilde{f}_{k}(\alpha) \ra^{2}_{\CH(K_1)} - 1 \nonumber\\
&=&  \left( \la \hat{\tilde{f}}_{kn}(\alpha) - \tilde{f}_{k}(\alpha) , \tilde{f}_{k}(\alpha) \ra_{\CH(K_1)} \right) \left( \la \hat{\tilde{f}}_{kn}(\alpha), \tilde{f}_{k}(\alpha) \ra_{\CH(K_1)} + 1 \right).
\end{eqnarray*}
Furthermore, we note that
\begin{eqnarray}\label{first identity in eigenvector proof} \sqrt{n} \left\{ \hat{\tilde{f}}_{kn}(\alpha) - \tilde{f}_{k}(\alpha) \right\} &=& \sqrt{n} \left[ \tilde{P}_{1k}(\alpha) \right] \left\{ \hat{\tilde{f}}_{kn}(\alpha) - \tilde{f}_{k}(\alpha) \right\}
\nonumber\\ &+& \sqrt{n} \left[I - \tilde{P}_{1k}(\alpha) \right] \left\{ \hat{\tilde{f}}_{kn}(\alpha) - \tilde{f}_{k}(\alpha) \right\}. \end{eqnarray}
Focussing on the first term in the right hand side of (\ref{first identity in eigenvector proof}) we see that
\begin{eqnarray}\label{first term expanded in eigenvector proof tikhinov section} &\sqrt{n}& \left[ \tilde{P}_{1k}(\alpha) \right] \left\{ \hat{\tilde{f}}_{kn}(\alpha) - \tilde{f}_{k}(\alpha) \right\}
\nonumber\\ &=& \left[ \sqrt{n} \la \hat{\tilde{f}}_{kn}(\alpha) - \tilde{f}_{k}(\alpha), \tilde{f}_{k}(\alpha) \ra_{\CH(K_1)} \right] \tilde{f}_{k}(\alpha) \nonumber\\
&=& \frac{  \sqrt{n} \left( \la \hat{\tilde{P}}_{1kn}(\alpha) - \tilde{P}_{1k}(\alpha), \tilde{P}_{1k}(\alpha) \ra_{HS} \right) }{\left( (\hat{\tilde{f}}_{kn}(\alpha), \tilde{f}_{k}(\alpha) )_{\CH(K_1)} + 1 \right)} \tilde{f}_{k}(\alpha). \end{eqnarray}
Due to the continuity of the inner product and (\ref{limiting distribution of P in Tikhinov section}) it follows that
\begin{eqnarray}\label{numerator in first term of tikhinov proof} &\sqrt{n}& \la \hat{\tilde{P}}_{1kn}(\alpha) - \tilde{P}_{1k}(\alpha), \tilde{P}_{1k}(\alpha) \ra_{HS} \nonumber\\
&\cid& \la \tilde{P}_{1k}(\alpha)  \CG_1(\alpha) \tilde{Q}_{1k}(\alpha), \tilde{P}_{1k}(\alpha) \ra_{HS} + \la \tilde{Q}_{1k}(\alpha)  \CG_1(\alpha) \tilde{P}_{1k}(\alpha), \tilde{P}_{1k}(\alpha) \ra_{HS} \nonumber\\
&=& {\rm tr} \left( \tilde{Q}_{1k}(\alpha) \CG_1(\alpha) \tilde{P}_{1k}(\alpha) \right) + {\rm tr} \left( \tilde{P}_{1k}(\alpha)  \CG_1(\alpha) \tilde{Q}_{1k}(\alpha) \tilde{P}_{1k}(\alpha) \right) \nonumber\\
&=& {\rm tr} \left( \tilde{Q}_{1k}(\alpha) \CG_1(\alpha) \tilde{P}_{1k}(\alpha) \right) + 0 \nonumber\\
&=& {\rm tr} \left( \tilde{P}_{1k}(\alpha) \tilde{Q}_{1k}(\alpha) \CG_1(\alpha)  \right) \nonumber\\
&=& 0
\end{eqnarray}
because
\begin{equation}\label{QPequal0} \tilde{Q}_{1k}(\alpha) \tilde{P}_{1k}(\alpha) = \sum_{j \neq k} \frac{1}{\rho_{j}(\alpha) - \rho_{k}(\alpha)} \tilde{P}_{1j}(\alpha) \tilde{P}_{1k}(\alpha) = 0 = \tilde{P}_{1k}(\alpha) \tilde{Q}_{1k}(\alpha).\end{equation}
Consequently, the numerator in (\ref{first term expanded in eigenvector proof tikhinov section}) converges in probability to $0$
whereas the denominator $\left( \la \hat{\tilde{f}}_{kn}(\alpha), \tilde{f}_{k}(\alpha) \ra_{\CH(K_1)} + 1 \right) = (2 + \CO_P(n^{-1/2}))$. Slutsky's theorem then implies that $\sqrt{n} \left[ \tilde{P}_{1k}(\alpha) \right] \left\{ \hat{\tilde{f}}_{kn}(\alpha) - \tilde{f}_{k}(\alpha) \right\} \cid 0$ and hence $\la \hat{\tilde{f}}_{kn}(\alpha), \tilde{f}_{k}(\alpha) \ra_{\CH(K_1)} \cip 1$.

To address the second term on the right hand side of (\ref{first identity in eigenvector proof}) we observe that as a consequence of Slutsky's theorem
\begin{eqnarray}\label{second term expanded in eigenvector proof tikhinov section} &\sqrt{n}& \left[ I - \tilde{P}_{1k}(\alpha) \right] \left\{ \hat{\tilde{f}}_{kn}(\alpha) - \tilde{f}_{k}(\alpha) \right\} \nonumber\\
&=& \sqrt{n} \left[ I - \tilde{P}_{1k}(\alpha) \right] \hat{\tilde{f}}_{kn}(\alpha) \nonumber\\
&=& \frac{\sqrt{n} \left[ I - \tilde{P}_{1k}(\alpha) \right] \left[ \hat{\tilde{P}}_{1kn}(\alpha) \right] \tilde{f}_{k}(\alpha)}{\la \hat{\tilde{f}}_{kn}(\alpha), \tilde{f}_{k}(\alpha) \ra_{\CH(K_1)}} \nonumber\\
&=& \frac{\sqrt{n} \left[ I - \tilde{P}_{1k}(\alpha) \right] \left[ \hat{\tilde{P}}_{1kn}(\alpha) - \tilde{P}_{1k}(\alpha) \right] \tilde{f}_{k}(\alpha)}{\la \hat{\tilde{f}}_{kn}(\alpha), \tilde{f}_{k}(\alpha) \ra_{\CH(K_1)}} \nonumber\\
&\cid&  \left[ I - \tilde{P}_{1k}(\alpha) \right] \left[  \tilde{P}_{1k}(\alpha)  \CG_1(\alpha) \tilde{Q}_{1k}(\alpha)  + \tilde{Q}_{1k}(\alpha)  \CG_1(\alpha) \tilde{P}_{1k}(\alpha) \right] \tilde{f}_{k}(\alpha) \nonumber\\
&=&  \left[ I - \tilde{P}_{1k}(\alpha) \right] \left[  \tilde{Q}_{1k}(\alpha)  \CG_1(\alpha) \right] \tilde{f}_{k}(\alpha) \nonumber\\
&=&  \tilde{Q}_{1k}(\alpha)  \CG_1(\alpha) \tilde{f}_{k}(\alpha).\end{eqnarray}
Equation (\ref{second term expanded in eigenvector proof tikhinov section}) establishes
(\ref{limiting distribution of eigenfunctions}) which completes the proof. \qed

We may now derive the limiting distribution for $\sqrt{n} \left[ \hat{\rho}_{kn}(\alpha) - \rho_{k}(\alpha) \right]$,
where $\{\hat{\rho}_{kn}(\alpha), \rho_{k}(\alpha)\}$ denotes the $k^{th}$ distinct eigenvalue associated
with $\{\hat{\CS}_{1n}(\alpha),\CS_{1}(\alpha)\}$.  In the following result, if $\rho_{k}(\alpha)$ has geometric
multiplicity $d_k$ then $\sqrt{n} \left[ \hat{\rho}_{kn}(\alpha) - \rho_{k}(\alpha) \right]$ will be regarded
as a vector of dimension $d_k$.

\begin{theorem}\label{asymptotic convergence of tikhinov canonical correlation}
Assume that ${\rm E} \| X \|_{L^2(E)}^4 < \infty$ and the $k^{th}$ regularized canonical correlation,
$\rho_{k}(\alpha)$, has geometric multiplicity $d_k$. Then,
\begin{eqnarray*}\label{asymptotic conv of eigenvalues of CS tikhinonv} \sqrt{n} \left[ \hat{\rho}_{kn}(\alpha) - \rho_{k}(\alpha) \right] &=& \sqrt{n} \left[ \hat{\tilde{P}}_{1kn}(\alpha) \hat{\CS}_{1n}(\alpha) \hat{\tilde{P}}_{1kn}(\alpha) -  \rho_{k}(\alpha) \tilde{P}_{1k}(\alpha) \right]
\nonumber\\ &\cid& \tilde{P}_{1k}(\alpha) \CG_1(\alpha) \tilde{P}_{1k}(\alpha) \end{eqnarray*}
with $\CG_1(\alpha)$ the Gaussian random variable in (\ref{complex distribution tikhinov section}). Furthermore,
$\tilde{P}_{1k}(\alpha) \CG_1(\alpha) \tilde{P}_{1k}(\alpha)$ has dimension $d_k$ and,
in the special case that $d_k = 1$,
\begin{equation*}\label{asymptotic conv of eigenvalues of CS tikhinonv 1 dim spec case} \sqrt{n} \left( \hat{\rho}_{kn}(\alpha) - \rho_{k}(\alpha) \right) \cid N(0,\sigma_{kk}(\alpha)) \end{equation*}
where $N(0,\sigma_{kk}(\alpha))$ denotes a normal distribution with zero mean and variance
\begin{equation*}\label{variance associated with eigenvalue} \sigma_{kk}(\alpha) = {\rm E} \left[ \la \tilde{f}_{k}(\alpha), \CG_1(\alpha) \tilde{f}_{k}(\alpha) \ra^2_{\CH(K_1)} \right]. \end{equation*}
\end{theorem}

\noindent Proof: Let $\rho_{k}(\alpha)$ denote the $k^{th}$ distinct eigenvalue of $\CS_{1}(\alpha)$ and assume
that it has multiplicity $d_k$. As $\| \hat{\tilde{P}}_{1kn}(\alpha) - \tilde{P}_{1k}(\alpha) \| \cip 0$,
Theorem \ref{continuity of specturm} ensures that for $n$ large enough,
$\rank(\hat{\tilde{P}}_{1kn}(\alpha)) = \rank(\tilde{P}_{1k}(\alpha)) = d_k$ with probability tending to $1$. Now observe that
\begin{eqnarray*}\label{eigenvalue decomposition defn in tikhinonv proof} \sqrt{n} \left[ \hat{\rho}_{kn}(\alpha) - \rho_{k}(\alpha) \right]
&=& \sqrt{n} \left[ \hat{\tilde{P}}_{1kn}(\alpha) \hat{\CS}_{1n}(\alpha) \hat{\tilde{P}}_{1kn}(\alpha) -  \rho_{k}(\alpha) \tilde{P}_{1k}(\alpha) \right]
\nonumber\\ &=& \sqrt{n} \left[ \sum_{j=1}^{3} \hat{\CB}_{kj}(\alpha) \right] \end{eqnarray*}
where
\begin{align*}\label{CB123 definition Tikhinov} \hat{\CB}_{k1}(\alpha) &\equiv   [ \hat{\tilde{P}}_{1kn}(\alpha) - \tilde{P}_{1k}(\alpha) ] \hat{\CS}_{1n}(\alpha) \hat{\tilde{P}}_{1kn}(\alpha), \notag\\
\hat{\CB}_{k2}(\alpha) &\equiv  \tilde{P}_{1k}(\alpha) [ \hat{\CS}_{1n}(\alpha) - \CS_{1}(\alpha) ] \hat{\tilde{P}}_{1kn}(\alpha), \notag\\
\hat{\CB}_{k3}(\alpha) &\equiv  \tilde{P}_{1k}(\alpha) \CS_{1}(\alpha) [ \hat{\tilde{P}}_{1kn}(\alpha) - \tilde{P}_{1k}(\alpha)].
\end{align*}
Equations (\ref{QPequal0}), (\ref{limiting distribution of P in Tikhinov section}) and Slutsky's theorem then ensure that
\begin{eqnarray*}\label{CB1 definition tikhinov} &&\| \hat{\CB}_{k1}(\alpha) \|^2_{HS}  \cid  \| \left[ \tilde{Q}_{1k}(\alpha) \CG_1(\alpha) \tilde{P}_{1k}(\alpha) + \tilde{P}_{1k}(\alpha) \CG_1(\alpha) \tilde{Q}_{1k}(\alpha) \right] \CS_{1}(\alpha) \tilde{P}_{1k}(\alpha) \|^2_{HS}  \nonumber\\
&\leq& \| \tilde{Q}_{1k}(\alpha) \CG_1(\alpha) \CS_{1}(\alpha) \tilde{P}_{1k}(\alpha) \|^2_{HS} + \| \tilde{P}_{1k}(\alpha) \CG_1(\alpha) \tilde{Q}_{1k}(\alpha) \CS_{1}(\alpha) \tilde{P}_{1k}(\alpha) \|^2_{HS} \nonumber\\
&=& \| \tilde{Q}_{1k}(\alpha) \CG_1(\alpha) \CS_{1}(\alpha) \tilde{P}_{1k}(\alpha) \|^2_{HS} + \| \tilde{P}_{1k}(\alpha) \CG_1(\alpha) \left[ \tilde{Q}_{1k}(\alpha)\tilde{P}_{1k}(\alpha) \right] \CS_{1}(\alpha)  \|^2_{HS} \nonumber\\
&=& \tr\left(\tilde{P}_{1k}(\alpha) \CS_{1}(\alpha) \CG_1(\alpha) \tilde{Q}^2_{1k}(\alpha) \CG_1(\alpha) \CS_{1}(\alpha) \tilde{P}_{1k}(\alpha) \right) + 0 \nonumber\\
&=& \tr\left( \CS_{1}(\alpha) \left[ \tilde{P}_{1k}(\alpha) \tilde{Q}_{1k}(\alpha) \right] \CG_1(\alpha) \CG_1(\alpha) \left[ \tilde{Q}_{1k}(\alpha) \tilde{P}_{1k}(\alpha) \right] \CS_{1}(\alpha) \right) = 0.\end{eqnarray*}
Similarly,
\begin{eqnarray*}\label{CB3 definition tikhinov} &&\| \hat{\CB}_{k3}(\alpha) \|^2_{HS}  \cid  \| \tilde{P}_{1k}(\alpha) \CS_{1}(\alpha) \left[ \tilde{Q}_{1k}(\alpha) \CG_1(\alpha) \tilde{P}_{1k}(\alpha) + \tilde{P}_{1k}(\alpha) \CG_1(\alpha) \tilde{Q}_{1k}(\alpha) \right] \|^2_{HS}  \nonumber\\
&\leq& \| \CS_{1}(\alpha) \left[ \tilde{P}_{1k}(\alpha) \tilde{Q}_{1k}(\alpha) \right] \CG_1(\alpha) \tilde{P}_{1k}(\alpha) \|^2_{HS} + \| \tilde{P}_{1k}(\alpha) \CS_{1}(\alpha) \CG_1(\alpha) \tilde{Q}_{1k}(\alpha) \|^2_{HS} \nonumber\\
&=& 0 +  \tr\left(\tilde{Q}_{1k}(\alpha)  \CG_1(\alpha) \CS_{1}(\alpha)  \tilde{P}_{1k}(\alpha) \CS_{1}(\alpha) \CG_1(\alpha) \tilde{Q}_{1k}(\alpha) \right) \nonumber\\
&=& \tr\left( \CG_1(\alpha) \left[ \tilde{Q}_{1k}(\alpha) \tilde{P}_{1k}(\alpha) \right] \CS_{1}(\alpha)  \tilde{P}_{1k}(\alpha) \CS_{1}(\alpha) \left[\tilde{P}_{1k}(\alpha) \tilde{Q}_{1k}(\alpha)\right] \CG_1(\alpha) \right) =0.\end{eqnarray*}
Hence Corollary \ref{complex theorem in tikhinov section} and Slutsky's Theorem ensure that
$$  \sqrt{n} \hat{\CB}_{k2}(\alpha) \cid  \tilde{P}_{1k}(\alpha) \CG_1(\alpha) \tilde{P}_{1k}(\alpha)$$
which proves the first part of the theorem.

To see the validity of the second part of the theorem, assume that $d_k = 1$ and observe that
\begin{equation*}\label{eigenvalue decomposition in tikhinonv proof} \sqrt{n} \left( \hat{\rho}_{kn}(\alpha) - \rho_{k}(\alpha) \right) = \sqrt{n} \left\{ \sum_{j=1}^{3} \hat{\CC}_{kj}(\alpha) \right\} \end{equation*}
where
\begin{align*}\label{CC123 definition tikhinov} \hat{\CC}_{k1}(\alpha) &\equiv  \la [ \hat{\tilde{f}}_{kn}(\alpha) - \tilde{f}_{k}(\alpha)]  ,  \hat{\CS}_{1n}(\alpha) \hat{\tilde{f}}_{kn}(\alpha) \ra_{\CH(K_1)}, \notag\\
\hat{\CC}_{k2}(\alpha) &\equiv   \la \tilde{f}_{k}(\alpha)  , [ \hat{\CS}_{1n}(\alpha) - \CS_{1}(\alpha)] \hat{\tilde{f}}_{kn}(\alpha) \ra_{\CH(K_1)}, \notag\\
\hat{\CC}_{k3}(\alpha) &\equiv  \la \tilde{f}_{k}(\alpha)  , \CS_{1}(\alpha) [ \hat{\tilde{f}}_{kn}(\alpha) - \tilde{f}_{k}(\alpha) ] \ra_{\CH(K_1)}.
\end{align*}
Note that $\hat{\CC}_{k1}(\alpha) \cip 0$ and $\hat{\CC}_{k3}(\alpha) \cip 0$ as a consequence of equations (\ref{limiting distribution of eigenfunctions}), (\ref{QPequal0}), and Slutsky's theorem since
\begin{eqnarray*}\label{CC1 goes to zero tikhinov} \hat{\CC}_{k1}(\alpha) &\cid&  \la \tilde{Q}_{1k}(\alpha)  \CG_1(\alpha) \tilde{f}_{k}(\alpha)  ,  \CS_{1}(\alpha) \tilde{f}_{k}(\alpha) \ra_{\CH(K_1)} \nonumber\\ &=& \sum_{j \neq k} \frac{\rho_{k}(\alpha)}{(\rho_{j}(\alpha) - \rho_{k}(\alpha))} \la \tilde{P}_{1j}(\alpha) \CG_1(\alpha) \tilde{f}_{k}(\alpha)  ,  \tilde{f}_{k}(\alpha) \ra_{\CH(K_1)} = 0 \end{eqnarray*}
and,
\begin{eqnarray*}\label{CC3 goes to zero tikhinov} \hat{\CC}_{k3}(\alpha) &\cid&  \la \tilde{f}_{k}(\alpha)  ,  \CS_{1}(\alpha) \tilde{Q}_{1k}(\alpha)  \CG_1(\alpha) \tilde{f}_{k}(\alpha) \ra_{\CH(K_1)} \nonumber\\ &=& \sum_{j \neq k} \frac{\rho_{k}(\alpha)}{(\rho_{j}(\alpha) - \rho_{k}(\alpha))} \la \tilde{P}_{1j}(\alpha) \tilde{f}_{k}(\alpha) ,  \CG_1(\alpha) \tilde{f}_{k}(\alpha)  \ra_{\CH(K_1)} = 0. \end{eqnarray*}
Application of Theorem \ref{complex theorem in tikhinov section} and Slutsky's Theorem then ensures that $\hat{\CC}_{k2}(\alpha) \cid \la \tilde{f}_{k}(\alpha), \CG_1(\alpha) \tilde{f}_{k}(\alpha) \ra_{\CH(K_1)}.$ Which completes the proof since
$${\rm E} \left[ \la \tilde{f}_{k}(\alpha), \CG_1(\alpha) \tilde{f}_{k}(\alpha) \ra_{\CH(K_1)} \right] = 0$$
and hence
$$\var\left[ \la \tilde{f}_{k}(\alpha), \CG_1(\alpha) \tilde{f}_{k}(\alpha) \ra_{\CH(K_1)} \right] = {\rm E} \left[ \la \tilde{f}_{k}(\alpha), \CG_1(\alpha) \tilde{f}_{k}(\alpha) \ra^2_{\CH(K_1)} \right] \equiv \sigma_{kk}(\alpha).$$
\qed

\noindent Now as the Gaussian operator $\CG_1(\alpha)$ is Hilbert-Schmidt, we note that the variances $\sigma_{kk}(\alpha) \downarrow 0$ as $k \uparrow \infty$.

Of natural interest is the degree of correlation between the regularized correlation estimators $\{\hat{\rho}_{kn}(\alpha), \hat{\rho}_{jn}(\alpha) \}$ with $j \neq k$. To investigate this association let us take the simple case where for $j \neq k$, the multiplicities $d_k = d_j = 1$. Then,
\begin{eqnarray*}\label{covariance between correlation estimators} &&\sigma_{jk}(\alpha) \equiv \cov[\hat{\rho}_{kn}(\alpha), \hat{\rho}_{jn}(\alpha) ]\nonumber\\ &=&
{\rm E} \left[ \la \hat{\tilde{f}}_{kn}(\alpha), \hat{\CS}_{1n}(\alpha) \hat{\tilde{f}}_{kn}(\alpha) \ra_{\CH(K_1)}  \la \hat{\tilde{f}}_{jn}(\alpha), \hat{\CS}_{1n}(\alpha) \hat{\tilde{f}}_{jn}(\alpha) \ra_{\CH(K_1)}  \right]
\nonumber\\ &=&
{\rm E} \left[ \la (\hat{\tilde{f}}_{kn}(\alpha) \otimes_{\CH(K_1)} \hat{\tilde{f}}_{kn}(\alpha) ), \left[ \hat{\CS}_{1n}(\alpha) \otimes_{HS} \hat{\CS}_{1n}(\alpha) \right] (\hat{\tilde{f}}_{jn}(\alpha) \otimes_{\CH(K_1)} \hat{\tilde{f}}_{jn}(\alpha)) \ra_{HS} \right]
\nonumber\\ &=&
\la (\tilde{f}_{k}(\alpha) \otimes_{\CH(K_1)} \tilde{f}_{k}(\alpha)), {\rm E} \left[ \hat{\CS}_{1n}(\alpha) \otimes_{HS} \hat{\CS}_{1n}(\alpha) \right] (\tilde{f}_{j}(\alpha) \otimes_{\CH(K_1)} \tilde{f}_{j}(\alpha)) \ra_{HS}
\nonumber\\ &=&
\la (\tilde{f}_{k}(\alpha) \otimes_{\CH(K_1)} \tilde{f}_{k}(\alpha)), \Sigma_1(\alpha) (\tilde{f}_{j}(\alpha) \otimes_{\CH(K_1)} \tilde{f}_{j}(\alpha) ) \ra_{HS_1} \equiv [\Sigma_1(\alpha)]_{jk}\end{eqnarray*}
where $\Sigma_1(\alpha) \equiv {\rm E} \left[ \hat{\CS}_{1n}(\alpha) \otimes_{HS} \hat{\CS}_{1n}(\alpha) \right],$
and $[\Sigma_1(\alpha)]_{jk}$ is the $\{j,k\}^{th}$ element of $\Sigma_1(\alpha)$.
These developments suggest that the $j^{th}$ and $k^{th}$ regularized canonical correlation estimators are not necessarily independent.

There are many similarities between the Tikhinov regularized version of canonical correlation analysis discussed here and those discussed in Cupidon et al. \citep{Cupidon06}.  However, it is important to distinctions between Cupidon et al. \citep{Cupidon06} and the method discussed here. Firstly, in Cupidon et al \citep{Cupidon06} the regularized operators discussed were of the form $(S_1 + \alpha I)^{-1/2} S_{12} (S_2 + \alpha I)^{-1} S_{21} (S_1 + \alpha I)^{-1/2}$ whereas in our approach they are $S_{12} (S_2 + \alpha I)^{-1} S_{21} (S_1 + \alpha I)^{-1}$.  Secondly, in Cupidon et al. \citep{Cupidon06} the operator approaches $RR^{*}$ as the regularization parameter approaches zero, and in our approach it tends to $TT^{*}$, an RKHS based operator which has well posed solutions on a closed domain.  Finally, the variance in the asymptotic distribution of $\hat{S}_{12} (\hat{S}_2 + \alpha I)^{-1} \hat{S}_{21} (\hat{S}_1 + \alpha I)^{-1}$ is the sum of 4 terms whereas the variance of $(\hat{S}_1 + \alpha I)^{-1/2} \hat{S}_{12} (\hat{S}_2 + \alpha I)^{-1} \hat{S}_{21} (\hat{S}_1 + \alpha I)^{-1/2}$ involves 5 terms.

\section{TSVD Regularization Approach}

In the Tikhinov approach to regularization, the operators $\{S_1, S_2\}$ are replaced with the operators $\{(S_{1}+\alpha I), (S_{2}+\alpha I)\}$ to obtain invertible operators. By contrast, the truncated singular value decomposition (TSVD) method of regularization replaces the compact operators $\{S_1, S_2\}$ with the finite rank (and rank-deficient) operators
\begin{eqnarray*}\label{TSVD Regularization with alpha} S_1(\alpha) &=&  \sum_{\lambda_{1i} > \alpha} \lambda_{1i} \phi_{1i} \otimes_{L^2(E_1)} \phi_{1i} \text{ and }
\nonumber\\ S_2(\alpha) &=&  \sum_{\lambda_{2i} > \alpha} \lambda_{2i} \phi_{2i} \otimes_{L^2(E_1)} \phi_{2i}.\end{eqnarray*}

Let us now define $m_{1}(\alpha) = \{\text{\# of } \lambda_{1i} > \alpha \}$ with similar definition holding for $m_{2}(\alpha)$. To ensure the equal dimensionality of the truncated versions of $S_1$ and $S_2$, it is advantageous to re-parameterize TSVD regularization in terms of $\{m_1(\alpha),m_2(\alpha)\}$, rather than $\alpha$.
In this regard, for simplicity we will always take $m = m_{1}(\alpha) = m_{2}(\alpha)$.
Notice that, under this re-parametrization,
the compact operators $S_1$ and $S_2$ are replaced by the finite dimensional operators $S_1(m) \equiv S_1 \Pi_1(m)$ and $S_2(m) \equiv S_2 \Pi_2(m)$, where $\Pi_1(m) \equiv \sum_{i=1}^m P_{1i}$ and $\Pi_2(m) \equiv \sum_{i=1}^m P_{2i}$, are the projection operators associated with the largest $m$ eigenvalues of $S_1$ and $S_2$ (or cumulative projection operators).
Much like $\alpha$ in Tikhinov regularization, the truncation parameter $m$ is the regularization parameter. In TSVD we are interested in the case where $m = m(\alpha) \rightarrow \infty$ which occurs when $\alpha \downarrow 0$. However, it is important to mention that TSVD regularization is widely used in statistical
practice, as it is common for a practicing statistician to discard right and left eigenvectors corresponding to small singular values after looking at, for example, a scree plot of the singular values.

Development of the theory associated with the TSVD version of regularized canonical correlation analysis can
now proceed along lines that are parallel to the developments in Sections 5 and 6.  Accordingly, let us define the operators $R(m): \CH_2 \mapsto \CH_1$ and $T(m): \CH(K_2) \mapsto \CH(K_1)$ by
\begin{equation}\label{Rm} R(m) \equiv S_1(m)^{1/2 \dag} S_{12} S_2(m)^{1/2 \dag} \end{equation}
and
\begin{equation}\label{Tm} T(m) \equiv \Gamma_1 S_1(m)^{1/2 \dag} S_{12} S_2(m)^{1/2 \dag} \Gamma_2^{-1}.\end{equation}
As all operators in (\ref{Rm}) and (\ref{Tm}) are bounded and $S_{12}$ is Hilbert-Schmidt,
both $R(m)$ and $T(m)$ are Hilbert-Schmidt. Also, as the regularization parameter $m \rightarrow \infty$, both $\Pi_1(m)$ and $\Pi_2(m)$ converge to the identity. Thus,
by Theorem \ref{HS Convergence theorem}, $R(m)$ converges in operator norm
to the operator $\Gamma^{-1}_1 T \Gamma_2$ and the continuity
of $\Gamma_1$ and $\Gamma_2$ then entail that
\beqs\label{Tk convergence} T(m) = \Gamma_1 R(m) \Gamma^{-1}_2 \rightarrow T ~~ \text{ as } ~~ m \rightarrow \infty \eeqs
with convergence in operator norm.

Now suppose that $\{\rho_j(m), f_{j}(m), g_{j}(m)\}_{j=1}^{\rank(R(m)^{*}R(m))}$ is the singular system for $R(m)$ with
\beqs\label{svd of Rk} R(m) = \sum_{i,j=1}^m  \frac{\gamma_{ij}}{\sqrt{\lambda_{1i}\lambda_{2j} }} \phi_{2j} \otimes_{\CH_2} \phi_{1i}
= \sum_{j=1}^{\rank(R(m)^{*}R(m))} \rho_j(m) \left[ g_j(m) \otimes_{\CH_2} f_j(m) \right]. \eeqs
Then,
\begin{eqnarray*}\label{svd of TK} T(m) = \Gamma_1 R(m) \Gamma^{-1}_2  &=& \sum_{j=1}^{\rank(R(m)^{*}R(m))} \rho_j(m) \left[ (\Gamma_2 g_j(m)) \otimes_{\CH(K_2)} (\Gamma_1 f_j(m)) \right]
\nonumber\\ &=& \sum_{j=1}^{\rank(T(m)^{*}T(m))} \rho_j(m) \left[ \tilde{g}_j(m)  \otimes_{\CH(K_2)} \tilde{f}_j(m) \right] \end{eqnarray*}
with $\{\tilde{f}_j(m)= \Gamma_1 f_j(m), \tilde{g}_j(m)=\Gamma_2 g_j(m)\}_{j=1}^{\rank(T(m)^{*}T(m))}.$ We note that in general $\rank(T(m)^{*}T(m)) \leq m$. Now by (\ref{RKHS Representation}) and because $\{\phi_{ij}\}_{j=1}^m$ are CONSs for $\left[ \Pi_i(m)\CH_i\right]$, for $i = 1,2$, it follows that the canonical weight functions in $\CH(K_1)$ and $\CH(K_2)$ may be written by
\begin{equation*}\label{TSVD canonical weight formula} \tilde{f}_j(m) = \sum_{i=1}^{m} \lambda_{1i} f_{ji}(m) \phi_{1i}  ~ \text{ and } ~
\tilde{g}_j(m) = \sum_{i=1}^{m} \lambda_{2i} g_{ji}(m) \phi_{2i} \end{equation*}
with
\begin{equation*}\label{TSVD canonical fourier coefficients} f_{ji}(m) = \la f_j(m),\phi_{1i} \ra_{\CH_1} ~ \text{ and } ~
g_{ji}(m) = \la g_j(m),\phi_{2j} \ra_{\CH_2}.\end{equation*}
The corresponding regularized canonical variables in $L^2_{X_1}$ and $L^2_{X_2}$ are
\begin{eqnarray*}\label{TSVD canonical variables formula} U_j(m) &= \Psi_1( \tilde{f}_j(m) ) = \sum_{i=1}^{m} f_{ji}(m) \la X_1,\phi_{1i} \ra_{\CH_2} \text{ and } \notag\\
V_j(m) &= \Psi_2( \tilde{g}_j(m) ) = \sum_{i=1}^{m} g_{ji}(m) \la X_2,\phi_{2i} \ra_{\CH_2}. \end{eqnarray*}

The TSVD parallel to Theorems \ref{CEGR THM1} and \ref{CEGR THM2} from Cupidon et al. \citep{Cupidon06} also hold.

\begin{theorem}\label{TSVD CEGR THM1}
For any $f \in \CH_1$ and $g \in \CH_2$, with $\tilde{f} = \Gamma_1(f)$, $\tilde{g} = \Gamma_2(g)$
\begin{equation}\label{TSVD limit2} 0 \leq \rho_j^2(m)  \uparrow \rho^2_j  \leq 1, \text{ as } m \rightarrow \infty.\end{equation}
\end{theorem}

\noindent Proof:
The convergence is from below since
$$\|T(m)\| = \|\Pi_1(m) T  \Pi_2(m) \|  < \|\Pi_1(m)\| ~ \| T \| ~ \|\Pi_2(m) \| < \|T\|.$$
To see that (\ref{TSVD limit2}) holds, fix $j \geq 1$ and observe that as $m \uparrow \infty$
\begin{equation*}\label{proof of tsvd eigenvalue convergence} \left( \rho^2_j - \rho_j^2(m) \right) = | \rho^2_j - \rho_j^2(m) |
\leq \| TT^{*} - T(m)T^{*}(m) \| \downarrow 0. \end{equation*} \qed

\begin{theorem}\label{TSVD CEGR THM2}
Let $\{\tilde{f}_j(m),\tilde{g}_j(m)\} = \{ \Gamma_1(f_j(m)),\Gamma_2(g_j(m))\} \in \{ \CH(K_1), \CH(K_2) \}$
denote the regularized weight functions corresponding to the TSVD version of canonical correlation analysis.
Then, as $m \rightarrow \infty$ for $j = 1,2,\ldots$,
\begin{equation*}\label{tsvd rspcc weight functions converge} \|\tilde{f}_j(m) - \tilde{f}_j \|_{\CH(K_1)} \rightarrow 0 ~ \text{ and }
 \| \tilde{g}_j(m) - \tilde{g}_j \|_{\CH(K_2)} \rightarrow 0.\end{equation*}
\end{theorem}

\noindent Proof: The proof here parallels the one for Theorem \ref{CEGR THM2}.  The idea is that
since $\|T(m) T^{*}(m) - T T^{*}\| \rightarrow 0$ as $m \rightarrow \infty,$
this implies that for any $j \in \BBZ$, the corresponding eigenprojection operators $\| P_j(m) - P_j \|_{HS} \rightarrow 0$.
If we now assume, WLOG, that $\la \tilde{f}_j(m),\tilde{f}_j \ra_{\CH(K_1)} \geq 0$, the relation
\begin{eqnarray*}\label{proof on convergence of eigenvectors}
\| P_j(m) - P_j \|^2_{HS} &=& \la P_j(m) - P_j, P_j(m) - P_j \ra_{HS}
\nonumber\\ &=& 2 - 2 \la P_j(m),P_j \ra_{HS}
\nonumber\\ &=& 2 - 2 \la (\tilde{f}_j(m) \otimes_{\CH(K_1)} \tilde{f}_j(m)),(\tilde{f}_j \otimes_{\CH(K_1)} \tilde{f}_j) \ra_{HS}
\nonumber\\ &=& 2 - 2 \la \tilde{f}_j(m),\tilde{f}_j \ra^2_{\CH(K_1)}
\nonumber\\ &=& 2(1- \la \tilde{f}_j(m),\tilde{f}_j \ra_{\CH(K_1)})(1+ \la \tilde{f}_j(m),\tilde{f}_j \ra_{\CH(K_1)})
\nonumber\\ &=& \|\tilde{f}_j(m) - \tilde{f}_j\|^2_{\CH(K_1)}(1+ \la \tilde{f}_j(m),\tilde{f}_j \ra_{\CH(K_1)})
\nonumber\\ &\geq& \|\tilde{f}_j(m) - \tilde{f}_j\|^2_{\CH(K_1)}
\end{eqnarray*}
implies that
\beqs\label{summary that eigenvectors converge tsvd proof}  \|\tilde{f}_j(m) - \tilde{f}_j\|^2_{\CH(K_1)} \leq \| P_j(m) - P_j \|_{HS} \rightarrow 0 \eeqs
as the regularization parameter $m \uparrow \infty$. \qed

Theorem \ref{TSVD CEGR THM2}, along with the continuity of the mappings $\Psi_1$ and $\Psi_2$,
ensures the convergence of the regularized canonical variables $U_j(m) = \Psi_1 ( \tilde{f}_j(m) )$ and $V_j(m) = \Psi_2 ( \tilde{g}_j(m))$
to the true canonical variables $\Psi_1(\tilde{f}_j)$ and $\Psi_2(\tilde{g}_j)$, as $m \rightarrow \infty$.

Let us now discuss the computation of the singular value
decomposition of $T(m)$.  To accomplish this it suffices to consider the eigenvalue-eigenvector decomposition of
$T(m)T^{*}(m)$. This is the finite rank operator given by
\begin{eqnarray}\label{TSVD exact operator for TTstar} \CT_1(m) &\equiv& T(m)T^{*}(m) = \Gamma_1 R(m) R^{*}(m) \Gamma^{-1}_1
\nonumber\\ &=& \Gamma_1 \Pi_1(m) S_1(m)^{1/2 \dag} S_{12} S_2(m)^{-1} S_{21} S_1(m)^{1/2 \dag} \Pi_1(m) \Gamma^{-1}_1. \end{eqnarray}
As was true for the Tikhinov case, problems arise from the presence of the unknown $\Gamma_1$ in $T(m)$.
However, unlike Tikhinov regularization the operators involved are finite rank.
Note that in the finite rank case $\Ima(S_1^{1/2}) = \overline{\Ima(S_1^{1/2})} = \ker(S_1)^{\perp}$ and we
may substitute $\Gamma_1 \Pi_1(m)$ with $S_1^{1/2}\Pi_1(m)$ directly.  Upon direct substitution of $S_1^{1/2}$ for $\Gamma_1$ in
(\ref{TSVD exact operator for TTstar}) we obtain the operator
\begin{equation}\label{TSVD after substitution for TTstar} \CS_1(m) =  \Pi_1(m) S_{12} S_2(m)^{\dag} S_{21} S_1(m)^{\dag} \end{equation}
which is a mapping from $\CH_1$ into $\CH_2$.
Much like its Tikhinov cousin, the operator $\CS_1(m)$ is self-adjoint since
\begin{eqnarray*}\label{TSVD CS1 self adjoint} \CS_1(m) &=& \sum_{i,j=1}^{m} \left(\frac{\gamma^2_{ij}}{\lambda_{1i} \lambda_{2j}} \right) \left[\phi_{1i} \otimes_{\CH_1} \phi_{1j} \right]
\nonumber\\ &=& S_1(m)^{\dag}  S_{12} S_2(m)^{\dag} S_{21} \Pi_1(m)  = \CS^{*}_1(m).\end{eqnarray*}
Additionally, since the operator is finite rank, it is Hilbert-Schmidt and admits the eigenvalue-eigenvector decomposition
\begin{equation*}\label{TSVD eigen expansion for CS}  \CS_1(m) = \sum_{j=1}^{\rank(\CS_1(m))} \rho^2_j(m) \left[ f_j(m) \otimes_{\CH_1} f_j(m) \right] \end{equation*}
with $\{\rho^2_j(m), f_j(m)\}_{j=1}^{\rank(\CS_1(m))},$ the eigensystem for $\CS_1(m)$.

All of the themes discussed in Section 2.1 are still applicable here.  For example, since the
eigenfunctions $\{f_j(m)\}$ are in $\Ima(S^{1/2})$, they belong to both $\CH(K_1)$ and $\CH_1$.  If the eigenfunctions are considered to be elements of $\CH(K_1)$ we will notate these as  $\{\tilde{f}_j(m)\}$ with $\{\tilde{f}_j(m) = f_j(m)\}.$ We will now show that as $m \uparrow \infty$
\begin{equation*}\label{TSVD proof after assignment} \Gamma_1^{-1} \CS_1(m) \Gamma_1 \rightarrow  \Gamma_1^{-1} T T^{*} \Gamma_1. \end{equation*}
To see this notice that
\begin{eqnarray*}\label{TSVD exact operator convergence for CS1} \Gamma_1^{-1} \CS_1(m) \Gamma_1  &\equiv& \sum_{i,j=1}^m \frac{\gamma^2_{ij}}{\lambda_{1i} \lambda_{2j}}  \left[ (\tilde{\phi}_{1i} \Gamma_1) \otimes_{\CH_1} (\Gamma^{-1}_1 \tilde{\phi}_{1j}) \right]
\notag\\ &\rightarrow& \sum_{i,j=1}^{\infty} \frac{\gamma^2_{ij}}{(\lambda_{1i})(\lambda_{2j})}
\left[ \phi_{1i} \otimes_{\CH_1} \phi_{1j} \right] = \Gamma_1^{-1} TT^{*} \Gamma_1. \end{eqnarray*}
Therefore, since $\CS_1(m)$ is Hilbert-Schmidt, Theorem \ref{HS Convergence theorem}
ensures that as the regularization parameter $m \rightarrow \infty$,
\begin{equation*}\label{TSVD cs convergence in norm} \| \Gamma_1^{-1} \CS_1(m) \Gamma_1 - \Gamma_1^{-1} TT^{*} \Gamma_1\|_{\CH_1}
= \| \CS_1(m)- TT^{*}\|_{\CH(K_1)} \rightarrow 0. \end{equation*}
If $\CS_1(m)$ is regarded as an operator on $\CH(K_1)$, we let $\{\rho^2_j(m),\tilde{f}_j(m)\}$ denote the eigenvalue and eigenvector pairs for the operator and using this notation
\begin{equation*}\label{TSVD eigen expansion for CS when considered RKHS based}  \CS_1(m) =
\sum_{j=1}^{\rank(\CS_1(m))} \rho^2_j(m) \left[ \tilde{f}_j(m) \otimes_{\CH(K_1)} \tilde{f}_j(m) \right]. \end{equation*}

\section{Asymptotics for the TSVD Operators}

In this section we will discuss the large sample distribution and consistency
of sample versions of $\CS_1(m)$.  The obvious estimator for this quantity is given by
\begin{equation}\label{CEGR sample estimator for TSVD operators} \hat{\CS}_{1n}(m) \equiv \hat{\Pi}_{1n}(m) \hat{S}_{12n} \hat{S}_{2n}(m)^{\dag} \hat{S}_{21n} \hat{S}_{1n}(m)^{\dag} \end{equation}
with $\hat{\Pi}_{in}(m) \equiv \sum_{j=1}^m \hat{P}_{ijn}$ and  $\hat{S}_{in}(m)^{\dag} = (\hat{S}_{in} \hat{\Pi}_{in}(m))^{\dag} = (\hat{S}_{in})^{\dag} \hat{\Pi}_{in}(m)$ for $i = 1,2$. By considering each factor associated with the TSVD operators in equation (\ref{CEGR sample estimator for TSVD operators}) it is clear that the asymptotic distribution will  differ from its corresponding Tikhinov counterpart. We begin our analysis by assuming that the joint process $X$ has zero mean and ${\rm E}[\|X\|^4_{L^2(E)} ] < \infty$. Accordingly,
we need to develop the asymptotic distribution of $\hat{\Pi}_{in}(m)$ and $\hat{S}_{in}(m)^{\dag}$ for $i = 1,2$.

\begin{corollary}\label{asymptotic distribution of pi}
Provided that ${\rm E}[ \|X\|^4_{L^2(E)} ] < \infty$, then for $i = 1,2$ and $m \geq 1$,
\begin{eqnarray}\label{asymp dist for Pi hat} &\sqrt{n}&(\hat{\Pi}_{in}(m) - \Pi_i(m)) \cid \CA_i(m) \equiv \sum_{j > m} \sum_{{k \neq j}\atop{k>m}} P_{ij} \CN_i P_{ik} + \sum_{k > m} \sum_{{j \neq k}\atop{j>m}} P_{ik} \CN_i P_{ij}
\nonumber\\ &=& \left[ I - \Pi_{i}(m) \right]\left(\sum_{j=1}^{\infty} (P_{ij} \CN_i Q_{ij} + Q_{ij} \CN_i P_{ij})  \right)\left[ I - \Pi_{i}(m) \right] \end{eqnarray}
where $\CN_i$ is a the distributional limit of $\sqrt{n} (\hat{S}_{in} - S_i).$
\end{corollary}

\noindent Proof: From Dauxois et al. \citep{DNR82} we know that for $i = 1,2$
\begin{equation}\label{limiting distribution of P} \sqrt{n} \left\{ \hat{P}_{ikn} - P_{ik} \right\} \cid  P_{ik}  \CN_i Q_{ik}  + Q_{ik}  \CN_i P_{ik} \end{equation}
where $\{\hat{P}_{ikn}, P_{ik}\}_{i=1}^{2},$ denote the eigenprojection operators for associated with the $k^{th}$ largest eigenvalues to $\{\hat{S}_i, S_i\}_{i=1}^2$
and
$$ Q_{ik} = \sum_{j \neq k} \frac{1}{\lambda_{ij} - \lambda_{ik} } P_{ij}.$$
For the asymptotic distribution of the cumulative eigenprojection operator  $\hat{\Pi}_{in}(m) = \sum_{j=1}^{m} \hat{P}_{ijn},$
we notice that for all $m \geq 1$ the cumulative sum of the first term on the right hand side of (\ref{limiting distribution of P}) is
$$ \sum_{j \leq m} P_{ij} \CN_i Q_{1j} = \sum_{j \leq k} P_{ij} \CN_i \sum_{k \neq j}  \frac{1}{\lambda_{ik} - \lambda_{ij} } P_{ik}$$
and involves terms like
\begin{equation*} \left[ \begin{array}{cccccc} 0 & \frac{P_{i1} \CN_i P_{i2}}{\lambda_{i2} - \lambda_{i1} }  & \frac{P_{i1} \CN_i P_{i3}}{\lambda_{i3} - \lambda_{i1} }  & \cdots & \frac{P_{i1} \CN_i P_{im}}{\lambda_{im} - \lambda_{i1} } & \cdots \\
\frac{P_{i2} \CN_i P_{i1}}{\lambda_{i1} - \lambda_{i2} } & 0 & \frac{P_{i2} \CN_i P_{i3}}{\lambda_{i3} - \lambda_{i2} }  & \cdots & \frac{P_{i2} \CN_i P_{1m}}{\lambda_{im} - \lambda_{i2} } & \cdots \\
\vdots & \vdots & \vdots & \ddots & \vdots & \cdots \\
\frac{P_{im} \CN_i P_{i1}}{\lambda_{i1} - \lambda_{im} } & \frac{P_{im} \CN_i P_{i2}}{\lambda_{i2} - \lambda_{im} }  & \frac{P_{im} \CN_i P_{13}}{\lambda_{i3} - \lambda_{im} }  & \cdots & 0 & \cdots \\
\end{array} \right]_{\cdot}  \end{equation*}
Hence for any term with $j,k \leq m$ and $j \neq k$ the upper triangular terms (UTT) involve $\frac{P_{ik} \CN_i P_{ij}}{\lambda_{ij} - \lambda_{ik} }$ and the lower triangular terms (LTT) are $ \frac{(P_{ij} \CN_1 P_{ik}) }{\lambda_{ik} - \lambda_{ij} } = \frac{ - (P_{ik} \CN_i P_{ij})^{*} }{\lambda_{ij} - \lambda_{ik} }$ so that $LTT = - UTT^{*}$.
Since
\begin{equation}\label{conv in dist for Pi hat proof 1}  \sqrt{n} (\hat{\Pi}_{in}(m) - \Pi_{i}(m)) \cid \sum_{j \leq m} P_{ij} \CN_i Q_{ij} +  \sum_{j \leq m} Q_{ij} \CN_i P_{ij} \end{equation}
and $(P_{ij} \CN_i Q_{ij})^{*} = Q_{ij} \CN_i P_{ij}$, the lower triangular terms in the first summand will cancel with the upper triangular terms in the second summand for all indices $i,j\leq m$. Equation (\ref{conv in dist for Pi hat proof 1}) then telescopes and produces the following new asymptotic result
\begin{eqnarray*}\label{conv in dist for Pi hat in proof 2} &\sqrt{n}&(\hat{\Pi}_i(m) - \Pi_i(m)) \cid \sum_{j > m} \sum_{{k \neq j}\atop{k>m}} P_{ij} \CN_i P_{ik} + \sum_{j > m} \sum_{{k \neq j}\atop{k>m}} P_{ik} \CN_i P_{ij}
\notag\\ &=& \left[ I - \Pi_{i}(m) \right]\left(\sum_{j=1}^{\infty} (P_{ij} \CN_i Q_{ij} + Q_{ij} \CN_i P_{ij})  \right)\left[ I - \Pi_{i}(m) \right]. \end{eqnarray*}\qed

\noindent It is important to note that Corollary \ref{asymptotic distribution of pi} has applications not just to canonical correlation analysis but also to principal component analysis.

Now, consider the asymptotic distribution of $\hat{S}_i(m)^{\dag} = (\hat{S}_{in} \hat{\Pi}_i(m))^{\dag}$ for some $m >0$ and $i = 1,2$.
In this regard, observe that the function $F(z) = z^{-1}$ is analytic everywhere in the complex plane except for a pole at zero. Therefore $F$ is analytic on the subset of the complex plane defined by $D_i \equiv \{ z \in \BBC: {\rm Re}(z) \geq  \lambda_{im} - \epsilon \}$ with $0< \epsilon < \lambda_{im}$. The set $D_i$ also
contains the spectrum of $S_i(m)$. Consequently, by the delta theorem (see Cupidon et al. \cite{Cupidon07} and appendix) we have that
\begin{equation}\label{convergence of Sim} \sqrt{n} \left\{\hat{S}_{in}(m)^{\dag} - S_i(m)^{\dag} \right\} \cid \CB_i(m) \end{equation}
for $i = 1,2$ where
\begin{equation}\label{CB1 and 2} \CB_i(m) \equiv  - \sum_{j=1}^{m} \lambda_{ij}^{-2} P_{ij} \CN_i P_{ij}
+ \sum_{{k \neq j} \atop {j,k \leq m}} \frac{\lambda_{ik}^{-1} - \lambda_{ij}^{-1}}{\lambda_{ik} - \lambda_{ij}} P_{ik} \CN_i P_{ij}. \end{equation}

The asymptotic analysis of $\sqrt{n}(\hat{\CS}_{in}(m) - \CS_{i}(m))$ for $i = 1,2$ may now proceed
where the application of the delta method leads to a product rule development.
For this purpose we introduce the following Gaussian Hilbert-Schmidt operators
\begin{eqnarray}\label{CF1} \CF_{11}(m) &\equiv& \CA_1(m) S_{12} S_2(m)^{\dag} S_{21} S_1(m)^{\dag}, \nonumber\\
\CF_{12}(m) &\equiv& \Pi_1(m) \CN_{12} S_2(m)^{\dag} S_{21} S_1(m)^{\dag}, \nonumber\\
\CF_{13}(m) &\equiv& \Pi_1(m) S_{12} \CB_2(m) S_{21} S_1(m)^{\dag},\nonumber\\
\CF_{14}(m) &\equiv& \Pi_1(m) S_{12} S_2(m)^{\dag} \CN_{21} S_1(m)^{\dag},\nonumber\\
\CF_{15}(m) &\equiv& \Pi_1(m) S_{12} S_2(m)^{\dag} S_{21} \CB_1(m),\nonumber\\
\CF_1(m) &\equiv& \sum_{j=1}^{5} \CF_{1j}(m) = \sum_{j=2}^{5} \CF_{1j}(m). \end{eqnarray}

\noindent The corollary below then results from the application of the delta theorem (see Cupidon et al. \cite{Cupidon07}).

\begin{corollary}\label{complex theorem in TSVD section}
If ${\rm E}[ \|X\|^4_{L^2(E)} ] < \infty$, then as $n \rightarrow \infty$,
\beq\label{complex distribution TSVD section} \sqrt{n}(\hat{\CS}_{1}(m) - \CS_{1}(m)) \cid \CF_1(m).\eeq
\end{corollary}

\noindent Proof: The proof follows along lines of the one for Corollary \ref{complex theorem in tikhinov section}.
Specifically, we begin by defining the elements
\begin{eqnarray*}\label{1st differences in proof for S TSVD asymptotics} \hat{\CA}_{11}(m) &\equiv& \left[ \hat{\Pi}_{1n}(m) - \Pi_{1}(m) \right] \hat{S}_{12n} \hat{S}_{2n}(m)^{\dag} \hat{S}_{21n} \hat{S}_{1n}(m)^{\dag},\notag\\
\hat{\CA}_{12}(m) &\equiv&  \Pi_{1}(m) \left[ \hat{S}_{12n} - S_{12} \right] \hat{S}_{2n}(m)^{\dag} \hat{S}_{21n} \hat{S}_{1n}(m)^{\dag},\notag\\
\hat{\CA}_{13}(m) &\equiv& \Pi_{1}(m) S_{12} \left[ \hat{S}_{2n}(m)^{\dag} - S_{2}(m)^{\dag} \right] \hat{S}_{21n} \hat{S}_{1n}(m)^{\dag},\notag\\
\hat{\CA}_{14}(m) &\equiv& \Pi_{1}(m) S_{12} S_{2}(m)^{\dag} \left[ \hat{S}_{21n} - S_{21} \right] \hat{S}_{1n}(m)^{\dag},\notag\\
\hat{\CA}_{15}(m) &\equiv& \Pi_{1}(m) S_{12} S_{2}(m)^{\dag} S_{21} \left[ \hat{S}_{1n}(m)^{\dag} - S_{1}(m)^{\dag} \right]. \end{eqnarray*}
With this notation we may write
\begin{equation*}\label{expansion of difference of cs in TSVD proof} \sqrt{n}(\hat{\CS}_{1n}(m) - \CS_{1}(m)) =
\sqrt{n} \left[ \sum_{j=1}^{5} \hat{\CA}_{1j}(m) \right]. \end{equation*}
The application of (\ref{conv in dist tikhinov section}), (\ref{asymp dist for Pi hat}), (\ref{convergence of Sim}) and Slutsky's Theorem then ensure that
\begin{equation*} \sqrt{n} \left[ \sum_{j=1}^{5} \hat{\CA}_{1j}(m) \right] \cid \CF_1(m)\end{equation*}
since, for example, the term $\hat{\CA}_{11}(\alpha)$ consists of the factor $\sqrt{n} \left[ \hat{\Pi}_{1n}(m) - \Pi_{1}(m) \right] \cid \CA_{1}(m)$
right-multiplied by the factor
\begin{equation*} \hat{S}_{12n} \hat{S}_{2n}(m)^{\dag} \hat{S}_{21n} \hat{S}_{1n}(m)^{\dag} \cip S_{12} S_{2}(m)^{\dag} S_{21} S_{1}(m)^{\dag}.\end{equation*}

We can now show that $\|\CF_{11}(m)\| = 0$ with probability 1.  To see this note that $\CF_{11}(m)$ is self-adjoint
because it is the distributional limit of self-adjoint operators. Furthermore, as a consequence of Corollary \ref{asymptotic distribution of pi}
\begin{eqnarray*}\label{Proof that CF11 is 0}  \CF_{11}(m) &=& \CA_1(m) S_{12} S_2(m)^{\dag} S_{21} S_1(m)^{\dag} \nonumber\\
&=& \left[ I - \Pi_{1}(m) \right] \CA_1(m) S_{12} S_2(m)^{\dag} S_{21} S_1(m)^{\dag} \left[ \Pi_{1}(m) \right]
\nonumber\\
&=& \left[ I - \Pi_{1}(m) \right] \CF^{*}_{11}(m) \left[ \Pi_{1}(m) \right]
\nonumber\\
&=& \left[ I - \Pi_{1}(m) \right]  S_1(m)^{\dag} S_{12} S_2(m)^{\dag}  S_{21} \CA_1(m) \left[ \Pi_{1}(m) \right]
\nonumber\\
&=& \left[ I - \Pi_{1}(m) \right] \left[ \Pi_{1}(m) \right] S_1(m)^{\dag} S_{12} S_2(m)^{\dag}  S_{21} \CA_1(m) \left[ I - \Pi_{1}(m) \right] \left[ \Pi_{1}(m) \right]
\nonumber\\
&=& 0. \end{eqnarray*}
Thus, $\|\CF_{11}(m)\| = 0$ with probability 1 and $\CF_1(m) = \sum_{j=2}^{5} \CF_{1j}(m)$. This completes the proof. \qed

Corollary \ref{complex theorem in TSVD section} ensures that for all $m \geq 1$
\beqs\label{consistency statement of cs TSVD section} \| \hat{\CS}_{1n}(m) - \CS_{1}(m) \| = \CO_{P}(n^{-1/2}) \cip 0.\eeqs
Hence, $ \hat{\CS}_{1n}(m)$  is consistent for $\CS_{1}(m)$.  The triangle inequality reveals the association between errors which originate from having a sample estimator and using regularization to approximate the desired operator $T T^{*}$,
\begin{equation}\label{consistency in the tsvd section to tstart}  \| \hat{\CS}_{1n}(m) - T T^{*} \|
\leq \| \hat{\CS}_{1n}(m) - \CS_{1}(m) \| + \| \CS_{1}(m) - T T^{*} \|. \end{equation}
The first term on the right-hand side of (\ref{consistency in the tsvd section to tstart}) is a random error that originates from using a sample estimator of $\CS_{1}(m)$ and tends to zero as $n \rightarrow \infty$.  Meanwhile, the second term on the right hand side of
(\ref{consistency in the tsvd section to tstart}) is a deterministic error that arises from using a regularized approximation of $T T^{*}$ and will tend to zero as $m \uparrow \infty$.

Since the limiting distribution for $\sqrt{n}(\hat{\CS}_{1n}(m) - \CS_{1}(m))$ has been established, we may derive large sample asymptotics for $\{\hat{\rho}_{1jn}(m),\hat{\tilde{P}}_{1jn}(m),\hat{\tilde{f}}_{1jn}(m)\}$
where these quantities represent the $j^{th}$ eigenvalue, eigenprojection and eigenvector for $\hat{\CS}_{1n}(m)$. Let $\{\rho_{1j}(m),\tilde{P}_{1j}(m),\tilde{f}_{1j}(m)\}$ denote similar quantities for $\CS_{1}(m)$.

We begin our development with the limiting distribution of the eigenprojection operators and associated eigenvectors.

\begin{theorem}\label{theorem for projections and eigenvectors TSVD section}
Suppose ${\rm E} \| X \|_{L^2(E)}^4 < \infty$. Then, for $m \geq 1$ and as $n \rightarrow \infty$
\beq\label{limiting distribution of P in TSVD section} \sqrt{n} \left\{ \hat{\tilde{P}}_{1kn}(m) - \tilde{P}_{1k}(m) \right\}
\cid  \tilde{P}_{1k}(m)  \CF_1(m) \tilde{Q}_{1k}(m)  + \tilde{Q}_{1k}(m)  \CF_1(m) \tilde{P}_{1k}(m) \eeq
where $\CF_1(m)$ is as in (\ref{complex distribution TSVD section}) and
\beqs\label{definition of Qj TSVD version}  \tilde{Q}_{1k}(m) = \sum_{j \neq k} \frac{1}{\rho_{1j}(m) - \rho_{1k}(m)} \tilde{P}_{1k}(m). \eeqs
In the case that $\rank(\tilde{P}_{1k}(m)) = 1$, then
\begin{equation}\label{limiting distribution of eigenfunctions TSVD version} \sqrt{n} \left\{ \hat{\tilde{f}}_{1kn}(m) - \tilde{f}_{1k}(m) \right\}
\cid  \tilde{Q}_{1k}(m)  \CF_1(m) \tilde{f}_{1k}(m). \end{equation}
\end{theorem}

\noindent Proof: The proof for the limiting distribution of $\hat{\tilde{P}}_{1kn}(m)$
is identical to that presented for Tikhinov regularization in Theorem \ref{theorem for projections and eigenvectors}. The only difference is that the role of the parameter $\alpha$ in Tikhinov regularization is replaced by that of $m$ in TSVD regularization. For the sake of completeness, we provide a sketch of the proof.

For each $k \geq 1$, let $\Gamma_k$ denote a circle that encloses the eigenvalue $\rho_{1k}(m)$ but no other eigenvalues of $\CS_{1}(m)$. From developments in the appendix, notice that
\begin{equation}\label{Asym 2 TSVD} \sqrt{n} \left\{ \hat{\tilde{P}}_{1kn}(m) - \tilde{P}_{1k}(m) \right\} =  \frac{\sqrt{n}}{2 \pi i} \oint_{\Gamma_j}  R(z)(\hat{\CS}_{1n}(m) - \CS_{1}(m)) R(z) dz + \CO_P(n^{-1/2}).\end{equation}
Focussing attention on the first term on the right hand side of (\ref{Asym 2 TSVD}), it follows from the continuous mapping theorem that
\begin{equation}\label{Asym 3 TSVD} \sqrt{n} \left\{ \hat{\tilde{P}}_{1kn}(m) - \tilde{P}_{1k}(m) \right\} \cid  \frac{1}{2 \pi i} \oint_{\Gamma_k}  R(z) \CF_1(m) R(z) dz. \end{equation}
Since
\begin{equation}\label{Resolvent formula} R(z) = \sum_{k=1}^{\infty} \frac{1}{\rho_{k}(m) - z} \tilde{P}_{1k}(m) + \CO( (\rho_{k}(m) -z)^{-2})\end{equation}
all but the lead term in (\ref{Resolvent formula}) will vanish when the contour integral is taken due to (\ref{Proof showing nilpotent terms can be ignored}).
The integrand in (\ref{Asym 3 TSVD}) can then be simplified as
\begin{equation*}\label{Sandwich} R(z) \CF_1(m) R(z)  = \sum_{k=1}^{\infty} \sum_{j=1}^{\infty}  \frac{1}{(\rho_{k}(m) - z)(\rho_{j}(m) - z)}\tilde{P}_{1k}(m) \CF_1(m) \tilde{P}_{1j}(m). \end{equation*}
Using the Cauchy integral formula produces
\begin{equation*}\label{integral} \frac{1}{2 \pi i} \oint_{\Gamma_k}  \frac{dz}{(\rho_{k}(m) - z)(\rho_{j}(m) - z)} = \frac{1}{2 \pi i} \oint_{\Gamma_k}  \frac{(\rho_{i}(m) - z)dz}{(\rho_{k}(m) - z)(\rho_{j}(m) - z)(\rho_{i}(m) - z)},\end{equation*}
and the only case where the integral is non-zero is when exactly one of $\rho_{k}(m)$ or $\rho_{j}(m)$ is not equal to $\rho_{i}(m)$. Hence
\begin{eqnarray*}\label{Asym Projections TSVD}  \sqrt{n} \left\{ \hat{\tilde{P}}_{1kn}(m) - \tilde{P}_{1k}(m) \right\} &\cid&
\sum_{i=1}^{\infty} \sum_{j \neq k} \frac{\delta_{ik}}{(\rho_{j}(m) - \rho_{k}(m))} \tilde{P}_{1i}(m) \CF_1(m) \tilde{P}_{1j}(m)  \nonumber\\ &+& \sum_{j=1}^{\infty} \sum_{i \neq k} \frac{\delta_{jk}}{(\rho_{i}(m) - \rho_{k}(m))} \tilde{P}_{1i}(m) \CF_1(m) \tilde{P}_{1j}(m)  \nonumber\\ &=& \tilde{P}_{1k}(m)  \CF_1(m) \tilde{Q}_{1k}(m)  + \tilde{Q}_{1k}(m)  \CF_1(m) \tilde{P}_{1k}(m)\end{eqnarray*}
which establishes (\ref{limiting distribution of P in TSVD section}).

To establish the limiting distribution of the eigenvectors in (\ref{limiting distribution of eigenfunctions TSVD version}) we write
\begin{eqnarray}\label{first identity in eigenvector proof tsvd case} \sqrt{n} \left\{ \hat{\tilde{f}}_{1kn}(m) - \tilde{f}_{1k}(m) \right\}
&=& \sqrt{n} \left[ \tilde{P}_{1k}(m) \right] \left\{ \hat{\tilde{f}}_{1kn}(m) - \tilde{f}_{1k}(m) \right\}
\nonumber\\ &+& \sqrt{n} \left[I - \tilde{P}_{1k}(m) \right] \left\{ \hat{\tilde{f}}_{1kn}(m) - \tilde{f}_{1k}(m) \right\}. \end{eqnarray}
Now by using the TSVD analogues to equations (\ref{first term expanded in eigenvector proof tikhinov section}) and (\ref{numerator in first term of tikhinov proof}) we may see that the limiting distribution for the first term on the right hand side of (\ref{first identity in eigenvector proof tsvd case}) is $0$. For the second term on the right hand side of (\ref{first identity in eigenvector proof tsvd case}) we have
$$ \sqrt{n} \left[ I - \tilde{P}_{1k}(m) \right] \left\{ \hat{\tilde{f}}_{1kn}(m) - \tilde{f}_{1k}(m) \right\} \cid \tilde{Q}_{1j}(m)  \CF_1(m) \tilde{f}_{1k}(m),$$
which completes the proof. \qed

We will now derive the limiting distribution for $\sqrt{n} \left[ \hat{\rho}_{1kn}(m) - \rho_{1k}(m) \right]$, where $\{\hat{\rho}_{1kn}(m), \rho_{1k}(m)\}$ denotes the $k^{th}$ distinct eigenvalue associated
with $\{\hat{\CS}_{1n}(m),\CS_{1}(m)\}$.  Much like the Tikhinov case, the quantity $\sqrt{n} \left[ \hat{\rho}_{1kn}(m) - \rho_{1k}(m) \right]$ will be regarded as a vector of dimension equal to the multiplicity, $d_k$, of the eigenvalue $\rho_{1k}(m)$.

\begin{theorem}\label{asymptotic convergence of TSVD canonical correlation}
Assume that ${\rm E} \| X \|_{L^2(E)}^4 < \infty$ and the $k^{th}$ regularized canonical correlation,
$\rho_{1k}(m)$, has geometric multiplicity $d_k$. Then
\begin{equation}\label{asymptotic conv of eigenvalues of CS TSVD section} \sqrt{n} \left[ \hat{\rho}_{1kn}(m) - \rho_{1k}(m) \right]
\cid \tilde{P}_{1k}(m) \CF_1(m) \tilde{P}_{1k}(m) \end{equation}
with $\CF_1(m)$ the Gaussian random variable in (\ref{CF1}). Furthermore,
$\tilde{P}_{1k}(m) \CF_1(m) \tilde{P}_{1k}(m)$ has dimension $d_k$.
In the special case that $d_k = 1$
\begin{equation}\label{asymptotic conv of eigenvalues of CS TSVD 1 dim spec case} \sqrt{n} \left( \hat{\rho}_{1kn}(m) - \rho_{1k}(m) \right) \cid N(0,\sigma_{kk}(m)) \end{equation}
where $N(0,\sigma_{kk}(m))$ denotes a normal distribution with zero mean and variance
\begin{equation*}\label{variance associated with eigenvalue TSVD section} \sigma_{kk}(m) = {\rm E} \left[ \la \tilde{f}_{1k}(m), \CF_1(m) \tilde{f}_{1k}(m)\ra^2_{\CH(K_1)} \right]. \end{equation*}
\end{theorem}

\noindent Proof: Like before, the proof here naturally parallels the Tikhinov result presented in Theorem \ref{asymptotic convergence of tikhinov canonical correlation}.
Since $\| \hat{\tilde{P}}_{1kn}(m) - \tilde{P}_{1k}(m) \| = \CO_P(n^{-1/2})$,
Theorem \ref{continuity of specturm} ensures that for large enough $n$,
$\rank(\hat{\tilde{P}}_{1kn}(m)) = \rank(\tilde{P}_{1k}(m)) = d_k$ with probability tending to $1$ as $n \rightarrow \infty$.
Now let us define
\begin{align*}\label{CD definition TSVD} \hat{\CD}_{k1}(m) &\equiv   [ \hat{\tilde{P}}_{1kn}(m) - \tilde{P}_{1k}(m) ] \hat{\CS}_{1n}(m) \hat{\tilde{P}}_{1kn}(m), \notag\\
\hat{\CD}_{k2}(m) &\equiv  \tilde{P}_{1k}(m) [ \hat{\CS}_{1n}(m) - \CS_{1}(m) ] \hat{\tilde{P}}_{1kn}(m), \notag\\
\hat{\CD}_{k3}(m) &\equiv  \tilde{P}_{1k}(m) \CS_{1}(m) [ \hat{\tilde{P}}_{1kn}(m) - \tilde{P}_{1k}(m)],
\end{align*}
and note that
\begin{equation}\label{rho decomposition} \sqrt{n} \left[ \hat{\rho}_{1kn}(m) - \rho_{1k}(m) \right] = \sqrt{n} \left[ \sum_{j=1}^{3} \hat{\CD}_{kj}(m)  \right].  \end{equation}
Note that $\hat{\CD}_{k1}(m) \cip 0$ and $\hat{\CD}_{k3}(m) \cip 0$ since
\begin{eqnarray*}\label{CD1 goes to zero} &&\| \hat{\CD}_{k1}(m) \|^2_{HS}  \cid  \| \left[ \tilde{Q}_{1k}(m) \CF_1(m) \tilde{P}_{1k}(m) + \tilde{P}_{1k}(m) \CF_1(m) \tilde{Q}_{1k}(m) \right] \CS_{1}(m) \tilde{P}_{1k}(m) \|^2_{HS}  \nonumber\\
&\leq& \| \tilde{Q}_{1k}(m) \CF_1(m) \CS_{1}(m) \tilde{P}_{1k}(m) \|^2_{HS} + \| \tilde{P}_{1k}(m) \CF_1(m) \tilde{Q}_{1k}(m) \CS_{1}(m) \tilde{P}_{1k}(m) \|^2_{HS} \nonumber\\
&=& \| \tilde{Q}_{1k}(m) \CF_1(m) \CS_{1}(m) \tilde{P}_{1k}(m) \|^2_{HS} + \| \tilde{P}_{1k}(m) \CF_1(m) \left[ \tilde{Q}_{1k}(m)\tilde{P}_{1k}(m) \right] \CS_{1}(m)  \|^2_{HS} \nonumber\\
&=& \tr\left(\tilde{P}_{1k}(m) \CS_{1}(m) \CF_1(m) \tilde{Q}^2_{1k}(m) \CF_1(m) \CS_{1}(m) \tilde{P}_{1k}(m) \right) + 0 \nonumber\\
&=& \tr\left( \CS_{1}(m) \left[ \tilde{P}_{1k}(m) \tilde{Q}_{1k}(m) \right] \CF_1(m) \CF_1(m) \left[ \tilde{Q}_{1k}(m) \tilde{P}_{1k}(m) \right] \CS_{1}(m) \right) = 0.\end{eqnarray*}
Similarly,
\begin{eqnarray*}\label{CD3 goes to zero} &&\| \hat{\CD}_{k3}(m) \|^2_{HS}  \cid  \| \tilde{P}_{1k}(m) \CS_{1}(m) \left[ \tilde{Q}_{1k}(m) \CF_1(m) \tilde{P}_{1k}(m) + \tilde{P}_{1k}(m) \CF_1(m) \tilde{Q}_{1k}(m) \right] \|^2_{HS}  \nonumber\\
&\leq& \| \CS_{1}(m) \left[ \tilde{P}_{1k}(m) \tilde{Q}_{1k}(m) \right] \CF_1(m) \tilde{P}_{1k}(m) \|^2_{HS} + \| \tilde{P}_{1k}(m) \CS_{1}(m) \CF_1(m) \tilde{Q}_{1k}(m) \|^2_{HS} \nonumber\\
&=& 0 +  \tr\left(\tilde{Q}_{1k}(m)  \CF_1(m) \CS_{1}(m)  \tilde{P}_{1k}(m) \CS_{1}(m) \CF_1(m) \tilde{Q}_{1k}(m) \right) \nonumber\\
&=& \tr\left( \CF_1(m) \left[ \tilde{Q}_{1k}(m) \tilde{P}_{1k}(m) \right] \CS_{1}(m)  \tilde{P}_{1k}(m) \CS_{1}(m) \left[\tilde{P}_{1k}(m) \tilde{Q}_{1k}(m)\right] \CF_1(m) \right) =0.\end{eqnarray*}
Hence, Corollary \ref{complex theorem in TSVD section} and Slutsky's Theorem ensure that
\begin{equation*}\label{CD2 convergence TSVD} \sqrt{n} \hat{\CD}_{k2}(m) \cid  \tilde{P}_{1k}(m) \CF_1(m) \tilde{P}_{1k}(m) \end{equation*}
which proves the first part of the theorem.

To see the validity of the second part of the theorem, assume that $d_k = 1$ and observe that
\begin{equation*}\label{eig decomp in tikhinonv proof expansion} \sqrt{n} \left( \hat{\rho}_{kn}(m) - \rho_{k}(m) \right) = \sqrt{n} \left\{ \sum_{j=1}^{3} \hat{\CC}_{kj}(m) \right\} \end{equation*}
where
\begin{align*}\label{CC123 TSVD} \hat{\CC}_{k1}(m) &\equiv  \la [ \hat{\tilde{f}}_{kn}(m) - \tilde{f}_{k}(m)]  ,  \hat{\CS}_{1n}(m) \hat{\tilde{f}}_{kn}(m) \ra_{\CH(K_1)}, \notag\\
\hat{\CC}_{k2}(m) &\equiv   \la \tilde{f}_{k}(m)  , [ \hat{\CS}_{1n}(m) - \CS_{1}(m)] \hat{\tilde{f}}_{kn}(m) \ra_{\CH(K_1)}, \notag\\
\hat{\CC}_{k3}(m) &\equiv  \la \tilde{f}_{k}(m)  , \CS_{1}(m) [ \hat{\tilde{f}}_{kn}(m) - \tilde{f}_{k}(m) ] \ra_{\CH(K_1)}.
\end{align*}
The terms $\hat{\CC}_{k1}(m) \cid 0$ and $\hat{\CC}_{k3}(m) \cid 0$ as a consequence of equation (\ref{limiting distribution of eigenfunctions TSVD version}) and Slutsky's
theorem, since
\begin{eqnarray*}\label{CC1 goes to zero TSVD} \hat{\CC}_{k1}(m) &\cid&  \la \tilde{Q}_{1k}(m)  \CF_1(m) \tilde{f}_{k}(m)  ,  \CS_{1}(m) \tilde{f}_{k}(m) \ra_{\CH(K_1)} \nonumber\\ &=& \sum_{j \neq k} \frac{\rho_{k}(m)}{(\rho_{j}(m) - \rho_{k}(m))} \la \tilde{P}_{1j}(m) \CF_1(m) \tilde{f}_{k}(m)  ,  \tilde{f}_{k}(m) \ra_{\CH(K_1)} = 0 \end{eqnarray*}
and
\begin{eqnarray*}\label{CC3 goes to zero TSVD} \hat{\CC}_{k3}(m) &\cid&  \la \tilde{f}_{k}(m)  ,  \CS_{1}(m) \tilde{Q}_{1k}(m)  \CF_1(m) \tilde{f}_{k}(m) \ra_{\CH(K_1)} \nonumber\\ &=& \sum_{j \neq k} \frac{\rho_{k}(m)}{(\rho_{j}(m) - \rho_{k}(m))} \la \tilde{P}_{1j}(m) \tilde{f}_{k}(m) ,  \CF_1(m) \tilde{f}_{k}(m)  \ra_{\CH(K_1)} = 0. \end{eqnarray*}

Application of Theorem \ref{complex theorem in TSVD section} implies that
\begin{equation*}\label{CD2 asym distr TSVD section} \sqrt{n} \hat{\CC}_{k2}(m) \cid \la \tilde{f}_{k}(m), \CF_1(m) \tilde{f}_{k}(m) \ra_{\CH(K_1)}.\end{equation*}
Since
\begin{equation*}\label{expectation of TSVD section} {\rm E} \left[ \la \tilde{f}_{k}(m), \CF_1(m) \tilde{f}_{k}(m) \ra_{\CH(K_1)} \right] = 0 \end{equation*}
and
\begin{equation*}\label{variance of TSVD section} \var\left[ \la \tilde{f}_{k}(m), \CF_1(m) \tilde{f}_{k}(m) \ra_{\CH(K_1)} \right] = {\rm E} \left[ \la \tilde{f}_{k}(m), \CF_1(m) \tilde{f}_{k}(m) \ra^2_{\CH(K_1)} \right] \equiv \sigma_{kk}(m),\end{equation*}
the proof is then complete. \qed

The TSVD versions of the correlation estimators $\{\hat{\rho}_{in}(m), \hat{\rho}_{jn}(m) \}$
with $i \neq j$ are correlated, much like the Tikhinov case. In fact, when the operator $\hat{\CS}_{1n}(m)$ is simple, we have for $i \neq j$
\begin{eqnarray*}\label{covariance between correlation estimators TSVD version} &&\sigma_{ij}(m) \equiv \cov[\hat{\rho}_{in}(m), \hat{\rho}_{jn}(m) ]
\nonumber\\ &=&
{\rm E} \left[ (\hat{\tilde{f}}_{in}(m), \hat{\CS}_{1n}(m) \hat{\tilde{f}}_{in}(m))_{\CH(K_1)}, (\hat{\tilde{f}}_{jn}(m), \hat{\CS}_{1n}(m) \hat{\tilde{f}}_{jn}(m))_{L^2(E_1)}  \right]
\nonumber\\ &=&
{\rm E} \left[ \left( \left(\hat{\tilde{f}}_{in}(m) \otimes_{\CH(K_1)} \hat{\tilde{f}}_{in}(m)\right), \left[ \hat{\CS}_{1n}(m) \otimes_{HS_1} \hat{\CS}_{1n}(m) \right] \left(\hat{\tilde{f}}_{jn}(m) \otimes_{\CH(K_1)} \hat{\tilde{f}}_{jn}(m) \right) \right)_{HS_1} \right]
\nonumber\\ &=&
\left( \left(\tilde{f}_{i}(m) \otimes_{\CH(K_1)} \tilde{f}_{i}(m)\right), {\rm E} \left[ \hat{\CS}_{1n}(m) \otimes_{HS_1} \hat{\CS}_{1n}(m) \right] \left(\tilde{f}_{j}(m) \otimes_{\CH(K_1)} \tilde{f}_{j}(m) \right) \right)_{HS_1}
\nonumber\\ &=&
\left( \left(\tilde{f}_{i}(m) \otimes_{\CH(K_1)} \tilde{f}_{i}(m)\right), \Sigma_1(m) \left(\tilde{f}_{j}(m) \otimes_{\CH(K_1)} \tilde{f}_{j}(m) \right) \right)_{HS_1}
\nonumber\\ &\equiv& [\Sigma_1(m)]_{ij}\end{eqnarray*}
with
\begin{equation*}\label{Sigma1m at the end} \Sigma_1(m) \equiv {\rm E} \left[ \hat{\CS}_{1n}(m) \otimes_{HS_1} \hat{\CS}_{1n}(m) \right]\end{equation*}
and $[\Sigma_1(m)]_{ij},$ the $\{i,j\}^{th}$ element of $\Sigma_1(m)$.

\section{Conclusion}

In Sections 5--9 we discussed how both Tikhinov and TSVD regularized estimators approach their intended target of approximation, the RKHS based operator $TT^{*}$, in the limits of their respective regularization parameters. We also showed that the asymptotics associated with Tikhinov and TSVD sample estimators $\{\hat{\CS}_{1n}(\alpha),\hat{\CS}_{1n}(m)\}$ are similar in the sense that for every distributional result for quantities relative to the Tikhinov estimator $\hat{\CS}_{1n}(\alpha)$, there is an analogous distributional result for its TSVD cousin $\hat{\CS}_{1n}(m)$.  The question to ask here is whether or not one form of regularization should be preferred over the other.

The answer to this question lies in one critical flaw in the Tikhinov approach to FCCA, which up to this point has not yet been discussed.  Although replacing the operators $\{S_1, S_2\}$ with $\{(S_1 + \alpha I), (S_2 + \alpha I)\}$ fixes the operators invertibility issues, the operators still theoretically have infinite dimensionality. Infinite dimensional operators are problematic because no computer will ever be able to estimate all the eigenvalues and eigenvectors. On the sample side, since the operators $\{\hat{S}_{1n}, \hat{S}_{2n}\}$ have rank at most $n$, they are rank deficient.  Meanwhile the operators $\{(\hat{S}_{1n} + \alpha I), (\hat{S}_{2n} + \alpha I)\}$ will have infinitely many eigenvalues equal to $\alpha$.  Any pragmatic computational scheme where Tikhinov regularization is implemented would therefore involve some limit on the number of eigenvalue and eigenvector pairs to be used and estimated.  As a consequence, FCCA methods will surely involve truncation. If we choose to implement Tikhinov regularization with truncation this will involve the operator
\begin{equation}\label{proposed estimator in conclusion} \hat{S}_{1n}(\alpha,m) \equiv  \sum_{j=1}^{m} (\hat{\lambda}_{1jn} + \alpha) \hat{P}_{1jn} = (\hat{S}_{1n} + \alpha I )\hat{\Pi}_{1n}(m) \end{equation}
for some integer $1 \leq m \leq n$.  The estimator in (\ref{proposed estimator in conclusion}) has some characteristics that are akin to both those of Tikhinov and TSVD regularization.  Utilizing this ``truncated Tikhinov'' estimator it follows that the corresponding regularized estimator for $TT^{*}$ would be
\begin{equation}\label{TSVD and Tikhinov based estimator for TTstar} \hat{\CS}_{1n}(\alpha,m) \equiv \hat{\Pi}_{1n}(m) \hat{S}_{12n} \hat{S}^{\dag}_{2n}(\alpha,m) \hat{S}_{21n} \hat{S}^{\dag}_{1n}(\alpha,m). \end{equation}
Equation (\ref{TSVD and Tikhinov based estimator for TTstar}) illustrates that pragmatic implementation of Tikhinov regularization in the FDA setting will in reality entail the use of both Tikhinov and TSVD forms of regularization. By contrast, TSVD regularization entails replacing the operators $\{S_1, S_2\}$ with $\{(S_1 \Pi_1(m) ,S_2 \Pi_2(m) \}$ which have finite rank. Consequently, TSVD regularization provides a remedy for both infinite dimensionality and invertibility issues simultaneously.

Since there are errors which originate from regularization methods in general, it is always better to use as few methods as possible.  The triangle inequality can now be utilized to establish a bound on the error associated with the ``truncated Tikhinov'' estimator (\ref{proposed estimator in conclusion}).  In this regard, notice that
\begin{equation*}\label{Triangle inequality in conclusion section}  \| \hat{\CS}_{1n}(\alpha,m) - T T^{*} \|
\leq  \| \hat{\CS}_{1n}(\alpha,m) - \hat{\CS}_{1n}(m) \| +  \| \hat{\CS}_{1n}(m) - T T^{*} \|. \end{equation*}
Hence the error associated with utilizing $\hat{S}_{1n}(\alpha,m)$ will always be larger than simply using $\hat{S}_{1n}(m)$.

\section{Appendix: Some Perturbation Theory}

In this appendix we briefly summarize some results from perturbation theory. The primary references for this section are Kato \citep{Kato80} and Dauxois et al. \citep{DNR82}.
A typical problem in perturbation theory is to determine how the eigenvalues and eigenspaces of a linear operator $B$ change when $B$ is subjected to a small perturbation.  Let $A: \CH \mapsto \CH$ be an arbitrary perturbation operator and let $\tilde{B} = B + A$ represent the perturbed operator.  In this regard, we might think of $A$ as being small in terms of its uniform operator norm $\|A\|$. However, a measure of ``closeness'' between $\tilde{B}$ and $B$ which is often of greater importance is the aperture or gap between the graphs of the two operators.

Let $\CM$ and $\CN$ be two closed linear manifolds on $\CH$ with $S_{\CM} = \{u \in \CM ~|~ \|u\|_{\CH} =1 \}$, the unit sphere on $\CM$. For any two closed linear manifolds $\CM, \CN \subset \CH$ let
\begin{equation*}\label{delta of manifolds} \delta (\CM, \CN) \equiv \left\{ \begin{array}{cc} {\sup}_{u \in S_{\CM}}\{\text{dist} (u, \CN)\} & \text{ for } \CM \neq \{0\}, \\ 0 & \text{ if } \CM = \{0\} \end{array} \right. \end{equation*}
with
\beqs\label{defn of distance} \text{dist} (u, \CN) \equiv \inf_{v \in \CN}\{ \| u - v\|_{\CH} \}.\eeqs
The gap between $\CM$ and $\CN$ is then defined by
\beqs\label{defn of gap} \hat{\delta}(\CM,\CN) \equiv \max[ \delta(\CM,\CN), \delta(\CN,\CM) ]. \eeqs
More details concerning $\delta(\CM,\CN)$ and $\hat{\delta}(\CM,\CN)$ can be found in Kato \citep{Kato80}.

If the graphs $\{G(B), G(\tilde{B})\}$ of two operators $\{B, \tilde{B}\}$ are closed, the closed graph theorem entails that both $B$ and $\tilde{B}$ are bounded. Consequently it is possible to define the gap between operators $B$ and $\tilde{B}$ by measuring the gap between their associated graphs.  In this regard we define
\beqs\label{defn of delta between ops} \delta(B,\tilde{B}) \equiv \delta(G(B), G(\tilde{B})), \eeqs
\beqs\label{defn of gap between ops} \hat{\delta}(B,\tilde{B}) \equiv \hat{\delta}(G(B), G(\tilde{B})) = \max[\delta(B,\tilde{B}),\delta(\tilde{B},B)], \eeqs
and $\hat{\delta}(B,\tilde{B}) = \hat{\delta}(\tilde{B},B)$ is called the gap between $B$ and $\tilde{B}$.

The notion of the gap between operators plays a large role in perturbation theory.  Suppose $B$ and $\tilde{B}$ are the original and perturbed operator respectively. The smaller the gap $\hat{\delta}(\tilde{B},B)$ becomes, the more properties the $\tilde{B}$ inherits from $B$.  Of particular importance is the following theorem from Kato \citep{Kato80} which permits the construction of closed curve $\Gamma$ around a part of the spectrum of $B$, denoted $\Sigma(B)$, that also encloses a similar collection of spectral points of the perturbed operator $\Sigma(\tilde{B})$.

\begin{theorem}\label{continuity of specturm}(Semi-continuity of the spectrum)
Let $\tilde{B}, B \in \CB(\CH)$ and let the spectrum of $B$, $\Sigma(B)$, be separated into two parts $\Sigma'(B)$, $\Sigma''(B)$ by a closed curve $\Gamma$, with $\CH = \CM'(B) \oplus \CM''(B)$. Then, there exists a $\delta > 0$, depending on $\Gamma$ and $B$, such that if $\tilde{B}$ is any operator with $\hat{\delta}(\tilde{B},B) < \delta$
\begin{enumerate}
\item the spectrum $\Sigma(\tilde{B})$ are likewise separated by $\Gamma$ into two parts $\{\Sigma'(\tilde{B}),\Sigma''(\tilde{B}) \}$ and both  $\{\Sigma'(\tilde{B}),\Sigma''(\tilde{B})\}$ are non-empty if this is true for $\{\Sigma'(B), \Sigma''(B)\}$,
\item in the associated decomposition $\CH = \CM'(\tilde{B}) \oplus \CM''(\tilde{B})$, $\{ \CM'(\tilde{B}), \CM''(\tilde{B})\}$ are isomorphic with $\{ \CM'(B), \CM''(B)\}$, respectively,
\item ${\rm{dim}}(\CM'(\tilde{B})) = {\rm{dim}}(\CM'(B))$ and ${\rm{dim}}(\CM''(\tilde{B})) = {\rm{dim}}(\CM''(B))$ \text{ and }
\item the projection operator $P_{\tilde{B}}$ of $\CH$ onto $\CM'(\tilde{B})$ tends to the similarly defined projection operator $P_{B}$ in operator norm as $\hat{\delta}(\tilde{B},B) \rightarrow 0$.
\end{enumerate}
\end{theorem}

We will now develop formulae for the differences in the resolvents and projection operators between the perturbed and unperturbed operator.  In this regard, let $R(z) = (B - z I)^{-1}$ and $\tilde{R}(z') = (\tilde{B} - z' I)^{-1}$ denote the resolvents of $B$ and $\tilde{B}$ for some $z \in \BBC \setminus \Sigma(B)$ and $z' \in \BBC \setminus \Sigma(\tilde{B})$, respectively.  From Kato \citep{Kato80}, if $\lambda_k \in \Sigma(B)$ is some isolated point of the spectra and $P_k$ is the associated projection operator then
\begin{equation}\label{Proj-res} P_k = - \frac{1}{2 \pi i} \oint_{\Gamma_k} R(z) dz  \end{equation}
where $\Gamma_k$ is a positively oriented curve that encloses  $\lambda_k$ but no other spectral values of $\Sigma(B)$.

Now, whenever $\|(\tilde{B}- B)R(z)\| < 1$ and $z \in \BBC \setminus \Sigma(B)$ we may utilize the Neumann series Theorem (Rynne and Youngson \cite{Rynne01}) which ensures that
\begin{eqnarray*}\label{Neuman series derivation of resolvant} \tilde{R}(z) &=& \left((\tilde{B}- B) + (B- zI) \right)^{-1} \nonumber\\ &=& \left((\tilde{B}- B) +(R(z))^{-1} \right)^{-1} \nonumber\\
&=& R(z) \left((\tilde{B}- B)R(z) + I \right)^{-1} \nonumber\\
&=& R(z) \left( I + \sum_{k=1}^{\infty} \left\{ (B - \tilde{B})R(z)\right\}^{k} \right). \end{eqnarray*}
It then follows that
\begin{eqnarray}\label{Difference between resolvants} \left[ \tilde{R}(z) - R(z) \right] &=& R(z)(B - \tilde{B}) R(z) \left[ \sum_{k=0}^{\infty} \left\{ (B - \tilde{B})R(z)\right\}^{k} \right]
\nonumber\\ &=&  R(z)(B - \tilde{B}) R(z) H(z)\end{eqnarray}
where $H(z) \equiv \sum_{k=0}^{\infty} \left\{ (B - \tilde{B})R(z)\right\}^{k}$.  Another application of Neumann series theorem reveals that
\begin{equation*}\label{H definition} H(z) = \sum_{k=0}^{\infty} \left\{ (B - \tilde{B})R(z)\right\}^{k} = \left[ I - (B - \tilde{B})R(z) \right]^{-1} \end{equation*}
and hence
\begin{equation*}\label{Difference between resolvants continued} \left[ \tilde{R}(z) - R(z) \right] =  R(z)(B - \tilde{B}) R(z) \left[ I - (B - \tilde{B})R(z) \right]^{-1}.\end{equation*}

Now let $\{\lambda_j,\tilde{\lambda}_j\}$ be particular spectral values for $\{B, \tilde{B}\}$, and let $\{P_j,\tilde{P}_j\}$ denote the corresponding eigenprojection operators. Now provided $\hat{\delta}(\tilde{B},B)$ is small enough, Theorem \ref{continuity of specturm} ensures that a positively oriented circle $\Gamma_j$, with radius $r$, can be drawn to enclose both  $\lambda_j$ and $\tilde{\lambda}_j$ but no other spectral values of either $B$ or $\tilde{B}$.  As a consequence of (\ref{Proj-res}) and (\ref{Difference between resolvants}) we then obtain
\begin{eqnarray}\label{Proj1} \left[ \tilde{P}_j - P_j \right] &=& - \frac{1}{2 \pi i} \oint_{\Gamma_j} \left[ \tilde{R}(z) - R(z) \right] dz \nonumber\\ &=& - \frac{1}{2 \pi i} \oint_{\Gamma_j}  R(z)(B - \tilde{B}) R(z) H(z) dz.\end{eqnarray}
Equation (\ref{Proj1}) allows us to formulate a crude bound on the uniform operator norm of  $\| \tilde{P}_j - P_j \|$, specifically
\begin{align}\label{Proj2} \| \tilde{P}_j - P_j \| & \leq  \frac{1}{2 \pi} \oint_{\Gamma_j} \| R(z)(B - \tilde{B}) R(z) \left[ I - (B - \tilde{B})R(z) \right]^{-1} \| dz
\notag\\ & \leq  \frac{2 \pi r}{2 \pi} \sup \left\{ \frac{\| B - \tilde{B} \| \| R(z)\|^2}{1 - \| B - \tilde{B} \| \| R(z)\| }:  z \in \Gamma_j \right\}
\notag\\ & =  r \sup \left\{ \frac{\| B - \tilde{B} \| \| R(z)\|^2}{1 - \| B - \tilde{B} \| \| R(z)\| }:  z \in \Gamma_j \right\} .\end{align}

Another formula for $\left[ \tilde{P}_j - P_j \right]$ can be derived by expanding the first term of $H(z)$ so that
\beq\label{H1} H(z) = I + \sum_{j=1}^{\infty} \left\{ (B - \tilde{B})R(z)\right\}^{j}. \eeq
Plugging (\ref{H1}) into (\ref{Proj1}) gives
\begin{equation}\label{Proj Diff} \left[ \tilde{P}_j - P_j \right] = - \frac{1}{2 \pi i} \oint_{\Gamma_j} R(z)(B - \tilde{B}) R(z) dz - \frac{1}{2 \pi i} \oint_{\Gamma_j} M(z) dz \end{equation}
where
\beqs\label{mn} M(z) \equiv R(z) \sum_{j = 2}^{\infty} \{ - A R(z) \}^{j} = \CO(A^2) .\eeqs
Using the partial fraction expansion of the resolvent from Kato \citep{Kato80} it follows that
\begin{equation}\label{R expansion} R(z) = \sum_{j=1}^{\infty} \frac{1}{(\lambda_j - z)} P_j + \CO( (\lambda_j -z)^{-2}). \end{equation}
Now in (\ref{R expansion}) the higher-order terms involving $\CO( (\lambda_j -z)^{-2})$ can be ignored due to Morera's theorem since, for $n \geq 2$,
\beq\label{Proof showing nilpotent terms can be ignored} \oint_{\Gamma_j} (\lambda_j -z)^{-n}  dz = \oint_{\Gamma_j} w^{n-2} dw = 0 \eeq
where the substitution $w = (\lambda_j - z)^{-1}$ has been used.
Thus, since
\beqs\label{cauchy1} \frac{1}{2 \pi i} \oint_{\Gamma_j} \frac{1}{\lambda_j - z} \frac{1}{\lambda_k - z} dz =
\left\{ \begin{array}{cc}   \frac{1}{\lambda_j - \lambda_k} & \text{ if } k \neq j, \\0 & \text{otherwise }\end{array} \right. \eeqs
it follows that
\begin{eqnarray*}\label{concl}  \frac{1}{2 \pi i} \oint_{\Gamma_j} R(z) A R(z) dz &=&  \frac{1}{2 \pi i} \sum_k \sum_j \oint_{\Gamma_j} \frac{1}{\lambda_k - z} \frac{1}{\lambda_j - z} dz P_k A P_j
\nonumber\\ &=& \sum_{k \neq l} \frac{1}{\lambda_j - \lambda_k} (P_k A P_j + P_j A P_k). \end{eqnarray*}
Therefore
\beq\label{big} \tilde{P}_j - P_j = \sum_{k \neq j} \frac{1}{\lambda_j - \lambda_k} (P_k A P_j + P_j A P_k) + \frac{1}{2 \pi i} \oint_{\Gamma} \phi(z) M(z) dz. \eeq

Equation (\ref{big}) has several important implications as it allows us to formulate the notion of the Frechet derivative of an analytic function of an operator. Now suppose that a function $\phi(z)$ is analytic in a domain $\Delta$ of the complex plane containing
all the spectral values $\{\lambda_h, \tilde{\lambda}_h\}$ of
$\{B ,\tilde{B}\}$, with $\Gamma \subset \Delta$ a positively
oriented closed curve that encloses all spectral values in its interior.  Utilizing the Dunsford-Taylor integral for $\phi(\tilde{B})$ and $\phi(B)$ (see Kato \citep{Kato80}) we obtain
\begin{eqnarray}\label{setup} \phi(\tilde{B}) - \phi(B) &=& - \frac{1}{2 \pi i} \oint_{\Gamma} \phi(z) [ \tilde{R}(z) - R(z) ] dz
\nonumber\\ &=& \frac{1}{2 \pi i} \oint_{\Gamma} \phi(z) R(z) A R(z) dz  + \frac{1}{2 \pi i} \oint_{\Gamma} \phi(z) M(z) dz
\nonumber\\ &=& \frac{1}{2 \pi i} \sum_k \sum_j \oint_{\Gamma}  \frac{\phi(z)}{(\lambda_k - z)(\lambda_j - z)} dz P_k A P_j + \CO(A^2). \end{eqnarray}
Focussing on the integral in the first term on the right hand side we see that
\begin{equation*}\label{first int} \frac{1}{2 \pi i} \oint_{\Gamma} \frac{\phi(z)}{(\lambda_k - z)(\lambda_j - z)} dz = \left\{ \begin{array}{cc}   \phi'(\lambda_j) & \text{ if } k = j, \\
\frac{\phi(\lambda_k)- \phi(\lambda_j)}{\lambda_k - \lambda_j} & \text{ if } k \neq j. \end{array}  \right. \end{equation*}
Equation (\ref{setup}) can then be written as
\begin{equation*}\label{setup2} \phi(\tilde{B}) - \phi(B) = \sum_{j \geq 1} \phi'(\lambda_j) P_j A P_j + \sum_{k \neq j}  \frac{\phi(\lambda_k)- \phi(\lambda_j)}{\lambda_k - \lambda_j}  P_k A P_j  + \CO(A^2). \end{equation*}
Now, since $\phi(\tilde{B}) = \phi(B) + \phi'_{B} A + \CO(A^2)$, the Frechet derivative at $B$ is
\beq\label{frechet deriv} \phi'_{B} A = \sum_{j \geq 1} \phi'(\lambda_j) P_j A P_j + \sum_{k \neq j}  \frac{\phi(\lambda_k)- \phi(\lambda_j)}{\lambda_k - \lambda_j}  P_k A P_j. \eeq
Equation (\ref{frechet deriv}) will be used extensively when we consider the delta method
for functions of random operators.

\bibliographystyle{apalike}

\end{document}